\documentclass[journal]{IEEEtran}
\IEEEoverridecommandlockouts                           

\usepackage{graphics} 
\usepackage{times} 
\usepackage{amsmath} 
\usepackage{amssymb}  
\usepackage{theorem}
\usepackage{color,xspace}
\usepackage{hyperref}
\usepackage{graphicx}
\usepackage{subcaption}
\usepackage{cite}
\usepackage{cleveref}
\usepackage{bm}
\usepackage{bbm}
\usepackage{tikz}
\usepackage{epstopdf} 

\input{mysymbol.sty}

\newtheorem{theorem}{Theorem}
\newtheorem{lemma}{Lemma}	

\newtheorem{remark}{Remark}

\newtheorem{proposition}{Proposition}	
\newtheorem{assumption}{Assumption}

\def\bu{\mathbf{u}}
\def\bx{\mathbf{x}}
\def\bz{\mathbf{z}}
\def\by{\mathbf{y}}

\def\bv{\mathbf{v}}

\def\bnu{\boldsymbol{\nu}}
\def\boldeta{\boldsymbol{\eta}}
\def\bbx{{\ensuremath{\mathbf x}}}

\title{Prediction-Correction Interior-Point Method for Time-Varying Convex Optimization}
\author{Mahyar Fazlyab, Santiago Paternain, Victor M. Preciado and Alejandro Ribeiro 
\thanks{This work was supported by the NSF under grants CNS-1302222, IIS-1447470 and by the ONR under grant N00014-12-1-0997. The authors are with the Department of Electrical and Systems Engineering, University of Pennsylvania. Email: \{mahyarfa, spater, preciado, aribeiro\}@seas.upenn.edu.}}

\begin{document}

\maketitle
\thispagestyle{empty}
\pagestyle{empty}

\begin{abstract}
%
	In this paper, we develop an interior-point method for solving a class of convex optimization problems with time-varying objective and constraint functions.
	Using log-barrier penalty functions, we propose a continuous-time dynamical system for tracking the (time-varying) optimal solution with an asymptotically vanishing error.
	This dynamical system is composed of two terms: (\emph{i}) a \emph{correction} term consisting of a continuous-time version of Newton's method, and (\emph{ii}) a \emph{prediction} term able to track the drift of the optimal solution by taking into account the time-varying nature of the objective and constraint functions.
	Using appropriately chosen time-varying slack and barrier parameters, we ensure that the solution to this dynamical system globally asymptotically converges to the optimal solution at an exponential rate.
	We illustrate the applicability of the proposed method in two practical applications: a sparsity promoting least squares problem and a collision-free robot navigation problem.
\end{abstract}



\section{Introduction}
The interplay between optimization and control theory is rich and fruitful, resulting in a plethora of efficient computational tools to solve fundamental control problems \cite{helmke2012optimization,feijer2010stability,wang2010control,wang2011control,gharesifard2014distributed,kia2015distributed,cherukuri2016asymptotic}.
Dynamical systems theory provides an array of mathematical tools to analyze the behavior of iterative algorithms proposed to solve standard optimization problems \cite{botsaris1978class,brown1989some}.
In this direction, control theory can be used to guarantee the convergence of iterative algorithms to accurate solutions, as well as to analyze the impact of numerical errors and computational delays.
Control tools have been extensively exploited in the context of stationary (i.e., time-invariant) optimization problems, in which both the objective function and constraints do not depend on time \cite{helmke2012optimization,feijer2010stability,wang2010control,wang2011control,gharesifard2014distributed,kia2015distributed,cherukuri2016asymptotic}.
In many practical settings, however, we find optimization problems in which the objective function and/or the constraints depend explicitly on time \cite{jakubiec2013d,cavalcante2013distributed,zhou2011multirobot,tu2011mobile, 7050337,su2009traffic,641452,716979}. In particular, time-varying optimization problems appear in, for example, the estimation of the path of a stochastic process \cite{jakubiec2013d}, signal detection with adaptive filters \cite{cavalcante2013distributed}, tracking of moving targets \cite{zhou2011multirobot}, and various problems in autonomous systems \cite{tu2011mobile, 7050337}, computer networks \cite{su2009traffic}, and learning \cite{641452,716979}.

In time-varying optimization problems, the optimal solution is a function of time; therefore, solving the optimization problem is equivalent to tracking the optimal solution as it varies over time.
A natural approach to addressing this problem is to sample the objective and constraint functions at particular times and to solve the corresponding sequence of (time-invariant) optimization problems using standard iterative algorithms (e.g., gradient or Newton's method \cite{boyd2004convex,fiacco1990nonlinear}).
However, this approach ignores the dynamic aspect of the problem, since each iteration tends to converge towards the optimal point of the sampled time-invariant problem, while the solution of the time-varying case is drifting away over time.
Therefore, this approach is likely to induce a steady-state optimality gap (i.e., a tracking error) whose magnitude depends on the time-varying aspects of the problem. This phenomenon has been previously observed in gradient descent algorithms for unconstrained optimization \cite{Popkov2005}, as well as in constrained optimization problems that arise in distributed robotics \cite{zavlanos2013network}, sequential estimation \cite{jakubiec2013d}, distributed optimization \cite{ling2013decentralized}, and neural networks \cite{zhang2011performance}.

In this paper, we consider time-varying smooth convex optimization problems characterized by (\emph{i}) a convex time-varying objective function, and (\emph{ii}) constraints that are expressed as level sets of time-varying convex functions and affine equalities. To track the time-varying optimal solution of the problem without a tracking error, we propose a continuous-time dynamical system whose state is globally asymptotically driven to the optimal solution at an exponential rate (under certain technical conditions).
In particular, we develop a prediction-correction interior-point method that utilizes information about time variations of the optimization problem in order to predict and correct the drift in the optimal solution, resulting in an asymptotically vanishing optimality gap.

The paper is structured as follows. In Section  \ref{se:statement}, we formally state the problem under consideration and introduce some regularity assumptions needed in our derivations. We then consider the particular case of time-varying optimization problems without constraints, as well as affine constraints only (Section \ref{sec_main_section}). In both cases, we propose to track the optimal solution using a dynamical system composed of two terms: (\emph{i}) a `prediction' term that uses information about time variations of the optimization problem, and (\emph{ii}) a `correction' term based on a continuous-time version of Newton's method. In Section \ref{sec_interior_point}, we propose a dynamical system able to track the solution of time-varying optimization problems with inequality constraints. In this case, we incorporate a logarithmic barrier function with an appropriately chosen time-varying barrier parameter, as well as a time-varying slack variable used to guarantee global convergence. We show that the proposed dynamical system converges at an exponential rate to the time-varying optimal point for any initial condition (Theorem \ref{thm_main_theorem}), under mild assumptions. These assumptions correspond to standard requirements to prove convergence of interior-point methods and differentiability of the objective and constraints with respect to time.
To illustrate our results, we perform numerical evaluations in a quadratic problem (Section \ref{sec_sim_sinthetic}) and consider two practical applications. The first application is a time-invariant $\ell_1$ regularized least squares problem in which we show that the use of a time-varying barrier parameter along with a prediction term speeds up the convergence of conventional interior-point methods (Section \ref{sec_static_interior_point}). The second application involves the navigation of a robot in an environment with circular obstacles (Section \ref{sec_robot_navigation}). We further consider situations in which the robot is charged with the task of tracking a moving target (Section \ref{sec_moving_target}). We close with concluding remarks (Section \ref{se:conclusions}).

\medskip\noindent{\it Notation.} Let $\mathbb{R}$, $\mathbb{R}_{+}$, and $\mathbb{R}_{++}$ be the set of real, nonnegative, and strictly positive numbers. The set $\{1,\ldots,n\}$ is denoted by $[n]$. We denote by $\mathbf{I}_n$ the $n$-dimensional identity matrix. We denote by $\mbS^{n}$ the space of $n$-by-$n$ symmetric matrices. The gradient of a function $f(\bx,t) \colon \mathbb{R}^n \times \mathbb{R}_{+} \to \mathbb{R}$ with respect to $\bx \in \mathbb{R}^n$ is denoted by $\nabla_{\bx} f(\bx,t)  \colon \mathbb{R}^n \times \mathbb{R}_{+} \to \mathbb{R}^n$. The \emph{partial} derivatives of $\nabla_{\bx} f(\bx,t)$ with respect to $\bx$ and $t$ are $\nabla_{\bx\bx} f(\bx,t)  \colon \mathbb{R}^n \times \mathbb{R}_{+} \to \mbS^n$ and $\nabla_{\bx t} f(\bx,t)  \colon \mathbb{R}^n \times \mathbb{R}_{+} \to \mathbb{R}^n$, respectively. \ifx The vector field $\mathbf{F}(\bx,t) \colon \mathbb{R}^n \times \mathbb{R}_{+} \to \mathbb{R}^n$ is said to be uniformly Lipschitz if there exists a positive constant $L$ such that $\|\mathbf{F}(\bx,t)-\mathbf{F}(\by,t)\|_2 \leq L\|\bx-\by\|_2$ for all $\bx,\by \in \mathbb{R}^n$ and all $t\in \mathbb{R}_{+}$.\fi


%
\section{Problem statement}\label{se:statement}

This paper considers a class of convex optimization programs where both the objective and the constraint functions are indexed by continuous time. Formally, consider a variable $\bbx\in\reals^n$ and let $t \in \mathbb{R}_{++}$ be a continuous time index. We then define a time-varying (TV) objective function $f_0:\reals^{n}\times\reals_+\to\reals$ taking values $f_0(\bbx,t)$; we also define $p$ time-varying inequality constraint functions $f_i:\reals^{n}\times\reals_+\to\reals$ taking values $f_i(\bbx,t)$ for $i\in [p]$; and $q$ time-varying  affine equality constraint functions $f'_i:\reals^{n}\times\reals_+\to\reals$ taking values $f'_i(\bbx,t)=\mathbf{a}_i(t)^\top \bx - b_i(t)$ for $i\in [q]$. For any given time $t>0$, the objective function and the constraints define an optimization problem whose optimal argument $\bbx^*(t)$ is defined as:
\begin{alignat}{2} \label{eq: inequality_constrained_time_varying_problem}
	\bbx^{\star}(t)
	:=\  &\argmin\ && f_{0}(\bbx,t),\\ 
	&\st    \ && f_{i}(\bbx,t) \leq 0, \quad i\in [p],   \nonumber \\
	& \ && \mathbf{A}(t) \bx = \mathbf{b}(t), \nonumber 
\end{alignat} 
where $\mathbf{A}(t) \colon \mathbb{R}_{+} \to \mathbb{R}^{q \times n}$ is defined as $\mathbf{A}^\top(t)=[\mathbf{a}_1(t) \cdots \mathbf{a}_q(t)]$,  and $\mathbf{b}(t)=[b_1(t) \cdots b_q(t)]^T$. A naive approach to solving \eqref{eq: inequality_constrained_time_varying_problem} is to sample the problem at particular times, say $0 \leq t_0<t_1<\cdots$, and solve the corresponding sequence of (time-invariant) optimization problems. In particular, for each $k \in \mathbb{Z}_{+}$, one could estimate $\bx^\star(t_k)$ by assuming the objective and constraints are time-invariant on the time interval $[t_k,t_{k+1})$ and run standard iterative algorithms (such as interior-point methods \cite{fiacco1990nonlinear,potra2000interior}). However, this naive implementation is likely to induce steady-state tracking errors. In particular, the sub-optimality of the solution computed at $t_{k+1}$ depends on the number of iterations allowed by our computational capabilities during the time interval $t \in [t_k,t_{k+1})$, as well as how fast the optimal argument has drifted away from $\bx^\star(t_k)$ during this interval.

Our goal is to develop an alternative approach to solving the TV optimization problem with a vanishing tracking error by leveraging information about the temporal variation of the objective and the constraints.
More precisely, we propose a continuous-time dynamical system $\dot{\bbx}(t) = \mathbf{h}(\bbx(t),t)$ whose solution $\bbx(t)$ satisfies $\|\bbx(t)-\bbx^{\star}(t)\| \to 0$ as $t\to \infty$ (i.e., it tracks the optimal solution with asymptotically vanishing error). To facilitate the exposition, we first address TV optimization problems without constraints (Section \ref{sec_unconstrained}), and then extend the framework to problems with TV constraints (Section \ref{sec_equality_constraint} and \ref{sec_interior_point}). In order to make the contributions of the paper more precise, we list below the assumptions that we impose on the optimization problem \eqref{eq: inequality_constrained_time_varying_problem}.

%
\begin{assumption}[Convexity] \label{assumption: convexity}
	\rm{The objective function $f_0(\bbx,t)$ and the constraint functions $f_i(\bbx,t),\ i\in [p]$ are twice continuously differentiable with respect to $\bbx$ and continuously differentiable with respect to $t$ for all $t \geq 0$. Furthermore, $f_i(\bbx,t)$ is convex in $\bbx$ for all $t \geq 0$ and $i\in [p]$.}
\end{assumption}
%
\begin{assumption}[Strong Convexity] \label{assumption: uniform_strong_convexity} \normalfont
	The objective function is uniformly strongly convex, i.e., there exists a positive constant $m$ such that for all times $t \geq 0$,  the inequality $\nabla_{\bbx\bbx} f_0(\bbx,t) \succeq m \mathbf{I}_n$ holds.
\end{assumption}
%
%
\begin{assumption}[Slater's condition] \label{assumption: strict_feasibility}
	\rm{The interior of the feasible region is nonempty for all $t \geq 0$, i.e., there exits $\bbx^\dagger \in \mbR^n$ such that $f_i(\bbx^\dagger,t)<0$, for all $i\in [p]$, and $\mathbf{a}_i(t)^\top \bx^{\dagger}=b_i(t),$ for all $i\in [q]$, $t \geq 0$.}
\end{assumption}
%
%
\begin{assumption} \label{assumption: equality constraints independent}
	\rm{The number of equality constraints is less than the dimension of the optimization space, i.e., $q<n$. Moreover, the vectors $\mathbf{a}_i(t),\ i\in [q]$ are linearly independent for all $t \geq 0$. This implies that $\mbox{rank}(\mathbf{A}(t))=q$ for all $t\geq 0$.}
\end{assumption}
%
%
%
%
\ifx
\green{
\begin{assumption} \label{assumption: uniform_lipschitz} \normalfont
The functions are uniformly Lipschitz ... I.e. there exists $L$ such that for all times $t$...
\end{assumption}}\fi

The uniform strong convexity of the objective function implies that the optimal trajectory $\bx^\star(t)$ is unique for all $t \geq 0$. From Assumption \ref{assumption: strict_feasibility}, the optimal solution $\bbx^\star(t)$ in \eqref{eq: inequality_constrained_time_varying_problem} at each $t \geq 0$ can be characterized using the Karush-Kuhn-Tucker (KKT) conditions \cite[Chapter 5]{boyd2004convex}. Finally, Assumption \ref{assumption: equality constraints independent} ensures that the system of equations $\mathbf{a}_i(t)^\top \bx = b_i(t)$ ($i\in [q]$) is consistent and has infinitely many solutions at each $t \geq 0$.
It is worth remarking that we do not make any assumption about asymptotic vanishing time variations in the objective and the constraints, namely, the partial derivates of these functions with respect to time are not assumed to converge to zero. However, we will assume that the optimal solution of the TV optimization problem does not grow exponentially as a function of time. We will explicitly impose this assumption in Section \ref{sec_interior_point}.


\section{Prediction-Correction Methods for Time-Varying Optimization}\label{sec_main_section}
In this section, we consider the TV optimization problem in \eqref{eq: inequality_constrained_time_varying_problem} without inequality constraints. In particular, we consider two versions of this problem: (\emph{i}) the unconstrained case (Section \ref{sec_unconstrained}), and (\emph{ii}) the case with linear equality constraints (Section \ref{sec_equality_constraint}). We show that, in both cases, it is possible to track the optimal trajectory $\bbx^{\star}(t)$ with an exponentially vanishing error. The algorithms developed here will be leveraged in Section \ref{sec_interior_point} to derive a prediction-correction interior-point method to tracking the solution of \eqref{eq: inequality_constrained_time_varying_problem} when inequality constraints are also considered.

%
\subsection{Unconstrained Time-Varying Convex optimization}\label{sec_unconstrained}

Consider the following unconstrained version of \eqref{eq: inequality_constrained_time_varying_problem}:
\begin{alignat}{2} \label{eq: unconstrained_time_varying_problem}
   \bbx^\star(t):= \argmin_{\bbx \in \reals^n} \ f_0(\bbx,t).  
\end{alignat} 
Under sufficient regularity conditions, we could implement a sequence of Newton steps on the function $f_0(\bbx,t)$  that would rapidly converge to $\bbx^\star(t)$. In the limit of infinitesimal steps, the sequence of Newton's iterations results in the following continuous-time dynamical system:
\begin{alignat}{2}\label{eqn_static_newton_dynamics}
\dot \bbx(t) = -\alpha \nabla_{\bbx\bbx}^{-1} f_0(\bbx(t),t) \nabla_{\bbx} f_0(\bbx(t),t),
\end{alignat}
where $\alpha >0$ is a constant. The trajectory $\bbx(t)$ generated by \eqref{eqn_static_newton_dynamics} would approach a neighborhood around $\bbx^\star(t)$, but it does {\it not} track exactly to $\bbx^\star(t)$, since the solution is itself changing over time. To overcome this limitation, observe that---with sufficient regularity---the optimal argument $\bbx^\star(t)$ in \eqref{eq: unconstrained_time_varying_problem} satisfies the first-order optimality condition $\nabla_{\bbx} f_0(\bbx^\star(t),t)=\bbzero$. Since this latter condition is true for all times $t \geq 0$, the time derivative of this condition should also be null, from where we obtain
\begin{alignat}{2}\label{eqn_optimum_following_dynamics_prelim}
   \bbzero  &= \dot \nabla_{\bbx} f_0(\bbx^\star(t),t)  \\ \nonumber
            &= \nabla_{\bbx\bbx} f_0(\bbx^\star(t),t) \dot\bbx^\star(t)
	              + \nabla_{\bbx t} f_0(\bbx^\star(t),t),
\end{alignat}
where $\dot \nabla_{\bbx} f_0$ denotes the total derivative of $\nabla_{\bbx} f_0$ with respect to time, while $\nabla_{\bbx t}f_0$ denotes the partial derivative of the gradient $\nabla_{\bbx} f_0$ with respect to time. Solving \eqref{eqn_optimum_following_dynamics_prelim} for $\dot\bbx^\star(t)$ yields the dynamical system
\begin{alignat}{2}\label{eqn_optimum_following_dynamics}
\dot\bbx^\star(t) = -\nabla_{\bbx\bbx}^{-1} f_0(\bbx^\star(t),t) \nabla_{\bbx t} f_0(\bbx^\star(t),t).
\end{alignat}
If the optimal solution $\bbx^\star(t)$ is known at some point in time, the system in \eqref{eqn_optimum_following_dynamics} can be used to track the evolution of $\bbx^\star(t)$, since \eqref{eqn_optimum_following_dynamics} guarantees that the optimality $\nabla_{\bbx} f_0(\bbx^\star(t),t)= \bbzero$ is satisfied for all times $t$. If we do not have access to $\bbx^\star(t)$ at any point in time, we propose to combine the dynamics in \eqref{eqn_static_newton_dynamics} and \eqref{eqn_optimum_following_dynamics} to build the following dynamical system:
\begin{align} \label{eq: time_varying_newton}
\dot \bbx(t) = -\nabla_{\bbx\bbx}^{-1} f_0(\bbx(t),t)
\Big[\alpha\nabla_{\bbx} f_0(\bbx(t),t) + \nabla_{\bbx t} f_0(\bbx(t),t) \Big].
\end{align}
The dynamics in \eqref{eq: time_varying_newton} contains two terms: (\emph{i}) a \emph{prediction} term $-\nabla_{\bbx\bbx}^{-1} f_0(\bbx(t),t)\nabla_{\bbx t} f_0(\bbx(t),t)$ that attempts to track the changes in the objective function [cf.  \eqref{eqn_optimum_following_dynamics}], and (\emph{ii}) a Newton-like \emph{correction} term $-\alpha\nabla_{\bbx\bbx}^{-1} f_0(\bbx(t),t)\nabla_{\bbx} f_0(\bbx(t),t)$ that `pushes' $\bbx(t)$ towards the optimum. In the next proposition, we show that the dynamical system in \eqref{eq: time_varying_newton} converges exponentially to the optimal trajectory $\bbx^\star(t)$.

%
\begin{proposition} \label{lem: Time_Varying_Newton}
	Let $\bx^{\star}(t)$ be defined as in \eqref{eq: unconstrained_time_varying_problem} and $\bx(t)$ be the solution of \eqref{eq: time_varying_newton}, where the objective function is $m$-strongly convex (Assumption \ref{assumption: uniform_strong_convexity}). Then, the following inequality holds:
\ifx
\begin{align}
\|\nabla_{\bx} f_0(\bx(t),t)\|_2 
\leq \|\nabla_{\bx} f_0(\bx(0),0)\|_2  e^{-\alpha t}.
\end{align} \fi
\begin{align}
	\|\bx(t)-\bx^{\star}(t)\|_2 
	     \leq \dfrac{1}{m}\|\nabla_{\bx} f_0(\bx(0),0)\|_2 e^{-\alpha t}.\nonumber
\end{align} 
\end{proposition}
%
%
\begin{myproof} See Appendix \ref{ap: time_varying_unconstrained_newton}.	\end{myproof}

%
Proposition \ref{lem: Time_Varying_Newton} confirms that the solution of \eqref{eq: time_varying_newton} tends exponentially to a point that satisfies the first-order optimality condition $\nabla_{\bbx} f_0(\bbx^\star(t),t)= \bbzero$.
For algorithmic implementation, we can discretize the continuous-time dynamics \eqref{eq: time_varying_newton} using, for example, Euler's forward method with a constant step size $\tau$ \cite[Chap. ~1]{iserles2009first}. If the vector field in the right-hand side of \eqref{eq: time_varying_newton} is uniformly Lipschitz in $\bx$, then the discretization error would be of the order $O(\tau)$ \cite[Chap. ~1]{iserles2009first}. Alternatively, one could use line search to generate adaptive step sizes that guarantee a strict reduction in the suboptimality in each iteration \cite{fazlyab2016self}.

\ifx
\blue{This is too vague: (i) What is Euler's method?, you can write it down or give a reference. What does convergence mean? Is the error in the solution of the differential equation? Or in the difference between the discretize solution and $\bbx^\star(t)$? Make it more explicit.} \blue{Also, this is perhaps better written as a remark that applies to all methods, instead of a remark about this particular method only.} \fi
%
\ifx
\begin{figure}\centering
\resizebox{9cm}{7cm}{
\input{figures/block_diagram_equality.tex}
}
\caption{Block diagram of the \red{write caption}} 
\label{fig_block_diagram_equality}
\end{figure}\fi
%
%
\subsection{Equality-Constrained Time-Varying Convex Optimization}\label{sec_equality_constraint}
We consider now a version of \eqref{eq: inequality_constrained_time_varying_problem} in which we incorporate equality constraints:
\begin{alignat}{2}\label{eq: equality_constraint_time_varying_problem}
  \bx^{\star}(t) := \argmin_{\bx \in \mathbb{R}^n} \  f_{0}(\bx,t),  \quad
                    \st\ \bbA(t)\bx=\mathbf{b}(t), 
\end{alignat} 
where the matrix $\bbA(t) \colon \mathbb{R}_{+} \to \mathbb{R}^{q\times n}$ and vector $\mathbf{b}(t) \colon \mathbb{R}_{+} \to \mathbb{R}^q$ define a TV equality constraint. 
In order to design a dynamical system to track $\bx^{\star}(t)$ in \eqref{eq: equality_constraint_time_varying_problem}, we introduce a Lagrange multiplier $\bnu \in \mathbb{R}^{q}$ and define the Lagrangian associated with the optimization problem in \eqref{eq: equality_constraint_time_varying_problem} as
\begin{align} \label{eq: time_varying_equality_lagrangian}
   \ccalL(\bx,\bnu,t):=f_0(\bx,t)+\bnu^\top(\bbA(t)\bx-\mathbf{b}(t)).
\end{align}
From the Lagrangian in \eqref{eq: time_varying_equality_lagrangian}, we define the TV dual function $\mathcal{G}(\bbnu,t):=\min_{\bbx\in\reals^n}\ccalL(\bx,\bnu,t)$ and the optimal dual argument as $\bbnu^*(t) := \argmax_{\bbnu\in\reals^q}\ \mathcal{G}(\bbnu,t)$.
By virtue of Assumptions \ref{assumption: uniform_strong_convexity} and \ref{assumption: equality constraints independent}, the optimal primal-dual pair $(\bx^{\star}(t),\bnu^{\star}(t))$ is unique for each time $t\geq 0$. Furthermore, we know this optimal pair must satisfy the following KKT conditions:
\begin{align}\label{eqn_equality_only_kkt}
   & \bbzero = \nabla_{\bx} \ccalL(\bx^{\star}(t),\bnu^{\star}(t),t)
             = \nabla_{\bx} f_0(\bx^\star(t),t)+\bbA(t)^\top \bnu^\star(t), \nonumber \\
   & \bbzero = \nabla_{\bnu} \ccalL(\bx^{\star}(t),\bnu^{\star}(t),t)
             = \bbA(t) \bx^\star(t)-\mathbf{b}(t) .
\end{align}
We define the aggregate variable $\bz :=[ \bx^\top, \bnu^\top]^\top \in \mathbb{R}^{n+q}$ and the optimal primal-dual solution $\bz^{\star}(t):=[ \bx^\star(t)^\top, \bnu^\star(t)^\top]^\top$ so as to rewrite \eqref{eqn_equality_only_kkt} in the condensed form~$\bbzero = \nabla_{\bz} \ccalL (\bx^{\star}(t),\bnu^\star(t),t) $. Since this latter equation must hold for all $t \geq 0$, we can take the time derivative of this gradient, which results in a prediction term of the form $-\nabla_{\bz\bz}^{-1} \ccalL(\bz(t),t)\nabla_{\bz t} \ccalL(\bz(t),t)$, similar to the one in \eqref{eqn_optimum_following_dynamics}. We combine this prediction term with a Newton-like correction term of the form $-\alpha \nabla_{\bz\bz}^{-1} \ccalL(\bz(t),t)\nabla_{\bz} \ccalL(\bz(t),t)$ to propose the following dynamical system:
\begin{align} \label{eq: time_varyin_equality_newton_method}
   \dot{\bz}(t)= -\nabla_{\bz\bz}^{-1} {\ccalL}(\bz(t),t)
                  \Big[\alpha \nabla_{\bz}{\ccalL}(\bz(t),t)
                         + \nabla_{\bz t}{\ccalL}(\bz(t),t)\Big].
\end{align}
In the following proposition, we prove that the state of this dynamical system converges exponentially to the optimal solution $\bx^{\star}(t)$ of \eqref{eq: equality_constraint_time_varying_problem}.
%
%
\begin{proposition} \label{lemma: equality_constrained_time_varyging_newton} Consider the optimization problem in \eqref{eq: equality_constraint_time_varying_problem} satisfying Assumptions
\ref{assumption: uniform_strong_convexity} and \ref{assumption: equality constraints independent}. Denote $\bz(t)=[\bx(t)^\top, \bnu(t)^\top]^\top$ as the solution of \eqref{eq: time_varyin_equality_newton_method},
where $ \ccalL(\bz,t)=\ccalL(\bx,\bnu,t)$ is defined in \eqref{eq: time_varying_equality_lagrangian}. Then, the following inequality holds:
\begin{align} \label{eq: time_varying_equality_constraint_gradient_exponential}
\| \nabla_{\bz} \ccalL(\bz(t),t)\|_2 \leq \| \nabla_{\bz} \ccalL(\bz(0),0)\|_2 e^{-\alpha t},
\end{align}
where $\bz(0)=[\bx(0)^\top \ \bnu(0)^\top]^\top \in \mathbb{R}^{n+p}$ is the initial condition  of \eqref{eq: time_varyin_equality_newton_method}. Furthermore, if the Hessian inverse of the Lagrangian is uniformly bounded (i.e., $\|\nabla_{\bz\bz}^{-1} L(\bz,t)\|_2 \leq M$ for some $M>0$ and all $t\geq 0$), the following holds true:
\begin{align} \label{eq: time_varying_equality_constraint_convergence_bound}
\|\bx(t)-\bx^{\star}(t)\|_2 &\leq Ce^{-\alpha t},
\end{align}
where $C=M\| \nabla_{\bz} \ccalL(\bz(0),0)\|_2$.
\end{proposition}
%
%
\begin{myproof} See Appendix \ref{ap: time_varying_equality_newton}. 
\end{myproof}
 
%
Proposition \ref{lemma: equality_constrained_time_varyging_newton} confirms that the state of the dynamical system in \eqref{eq: time_varyin_equality_newton_method} converges exponentially to the optimal solution $\bx^{\star}(t)$ of \eqref{eq: equality_constraint_time_varying_problem}. Notice also that this dynamical system can have any arbitrary initial condition; in other words, the state of the system is \emph{globally} asymptotically driven to the TV optimal solution  $\bx^{\star}(t)$ at an exponential rate. In fact, asymptotic feasibility is achieved, since the term $\|\nabla_{\bnu} \ccalL(\bx(t),\bnu(t),t)\|_2=\|A(t)\bx(t)-\mathbf{b}(t)\|_2$, which quantifies the infeasibility, vanishes exponentially, according to \eqref{eq: time_varying_equality_constraint_gradient_exponential}. Moreover, we show below that, if the dynamical system in \eqref{eq: time_varyin_equality_newton_method} starts from a feasible point $\bz(0)$, it will remain feasible for all $t>0$. To prove this, notice that the solution of the dynamical system in \eqref{eq: time_varyin_equality_newton_method} satisfies 
$\dot{\nabla}_{{\bnu}} \mathcal{L}(\bz(t),t)=-\alpha \nabla_{{\bnu}} \mathcal{L}(\bz(t),t)$, as we show in the proof of Proposition \ref{lemma: equality_constrained_time_varyging_newton} (see Eq. \eqref{eq: Lagrangian_gradient_dynamics} in  Appendix \ref{ap: time_varying_equality_newton}). Therefore, since $\nabla_{\bnu} \ccalL(\bx(t),\bnu(t),t)=\mathbf{A}(t) \bx(t)-\mathbf{b}(t)$, we obtain:
\begin{align}
\frac{d}{dt} (\mathbf{A}(t) \bx(t)-\mathbf{b}(t)) = -\alpha (\mathbf{A}(t)\bx(t)-\mathbf{b}(t)).\nonumber
\end{align}
This implies that, if the solution is initially feasible (i.e., $\mathbf{A}(0) \bx(0)-\mathbf{b}(0)=\mathbf{0}$), then $\mathbf{A}(t) \bx(t)-\mathbf{b}(t)=\mathbf{0}$ for all $t \geq 0$. 

\begin{remark} \normalfont 
	In contrast to the unconstrained problem considered in Subsection \ref{sec_unconstrained}, where the objective function $f_0(\bx,t)$ is strongly convex in the primal variable $\bx \in \mathbb{R}^n$, the Lagrangian function of the problem with equality constraints, i.e., $\ccalL(\bx,\bnu,t)$, is no longer strongly convex in the primal-dual variable $\bz=[\bx^\top \ \bnu^\top]^\top \in \mathbb{R}^{n+q}$.\end{remark}
	
\begin{remark} \normalfont 
The assumption $\|\nabla_{\bz\bz}^{-1} L(\bz,t)\|_2 \leq M$ in Proposition \ref{lemma: equality_constrained_time_varyging_newton} allows to quantify the convergence rate in the domain of the primal variable $\bx$ [cf. \eqref{eq: time_varying_equality_constraint_convergence_bound}]. The latter assumption is also of practical importance for the following reason: if $\dot{\bz}(t)$ in \eqref{eq: time_varyin_equality_newton_method} grew arbitrarily large, then the discretization of \eqref{eq: time_varyin_equality_newton_method} would require arbitrarily small step sizes.
\end{remark}

In the next subsection, we propose a solution to the most general time-varying optimization problem, where both equality and inequality constraints are considered.

\subsection{General Time-Varying Convex Optimization }\label{sec_interior_point} 
In this subsection, we return to our original optimization problem in \eqref{eq: inequality_constrained_time_varying_problem}, considering both equality and inequality constraints. In light of the analysis in Subsection \ref{sec_equality_constraint}, we can always eliminate the equality constraints in \eqref{eq: inequality_constrained_time_varying_problem} by Lagrangian relaxation. Therefore, we ignore equality constraints for now, without losing generality, and will remark on the addition of equality constraints at the end of this subsection. Under Assumptions \ref{assumption: convexity} and \ref{assumption: strict_feasibility}, the necessary and sufficient KKT conditions \cite[Chapter 5]{boyd2004convex} for optimality of $\bbx^\star(t)$ in \eqref{eq: inequality_constrained_time_varying_problem} at each $t \geq 0$ read as,
\begin{align} \label{eq: time_varying_KKT_conditions}
& \nabla_{\bbx} f_0(\bbx^{\star}(t),t)+\sum_{i=1}^{p}\lambda_{i}^{\star}(t) \nabla_{\bbx} f_{i}(\bbx^{\star}(t),t)=\bbzero; \nonumber \\	
& \lambda_{i}^{\star}(t) f_{i}(\bbx^{\star}(t),t)=0,\ \lambda_{i}^{\star}(t) \geq 0,  \nonumber \\
& f_{i}(\bbx^{\star}(t),t) \leq 0,\text{ for all } i \in [p].
\end{align}
%
In what follows, we use barrier functions \cite[Chapter 11]{boyd2004convex} to incorporate the inequality constraints into the objective function. First, we consider the following convex optimization problem, which is equivalent to \eqref{eq: inequality_constrained_time_varying_problem}, without equality constraints:
\begin{alignat}{2} \label{eq: time_varying_problem_with_barrier}
	\bbx^{\star}(t)
	:=\  &\underset{\bx \in \mathbb{R}^n}{\argmin}\ && f_{0}(\bbx,t)+\sum_{i=1}^{p} \mathbb{I}_{-}(f_{i}(\bbx,t)),
\end{alignat}
where $\mathbb{I}_{-}: \mathbb{R} \to \{0,\infty\}$ is defined such that $\mathbb{I}_{-}(u)=0$ if $u\leq 0$, and $\mathbb{I}_{-}(u)=\infty$ if $u>0$. We now approximate $\mathbb{I}_{-}(u)$ by a (smooth) barrier function of the form $-\frac{1}{c} \log(-u)$, where $c>0$ is an arbitrary constant called the barrier parameter\footnote{Notice that $\lim_{c \to \infty} -\frac{1}{c} \log(-u) = \mathbb{I}_{-}(u)$.}. Therefore, we can approximate \eqref{eq: time_varying_problem_with_barrier} by the following smooth convex optimization problem:
\begin{alignat}{2} \label{eq: time_varying_problem_with_log_barrier}
    \underset{\bx \in \mathcal{D}(t)}{\mbox{minimize}}\quad && f_{0}(\bbx,t)-\dfrac{1}{c(t)}\sum_{i=1}^{p} \log(-f_{i}(\bbx,t)),
\end{alignat}
where $c(t)$ is a time-dependent positive barrier parameter and the domain of the objective function is the open set $\mathcal{D}(t):= \{\bbx \in \mathbb{R}^n \colon f_{i}(\bbx,t)<0,\ i\in [p] \}$. 

Our goal is to design a dynamical system able to track the optimal solution of \eqref{eq: time_varying_problem_with_log_barrier}.
As we show below, this would require the initial condition of this dynamical system to lie inside the initial domain $\mathcal{D}(0)$, i.e., $\bbx(0)\in\mathcal{D}(0)$.
To circumvent this requirement, we include a nonnegative time-dependent slack variable $s(t)$ in the optimization problem \eqref{eq: time_varying_problem_with_log_barrier}, and solve for the \emph{approximate} optimal trajectory defined by
\begin{alignat}{2} \label{eq: time_varying_problem_with_log_barrier_and_slack}
\widehat{\bbx}^{\star}(t)
:=\  &\underset{\bx \in \widehat{\mathcal{D}}(t)}{\argmin}\ && f_{0}(\bbx,t)-\dfrac{1}{c(t)}\sum_{i=1}^{p} \log(s(t)-f_{i}(\bbx,t)),
\end{alignat}
where $\widehat{\mathcal{D}}(t):= \{\bbx \in \mathbb{R}^n \colon f_{i}(\bbx,t)<s(t),\ i\in [p] \}$. Notice that for any $\bx(0) \in \mathbb{R}^n$, we can choose $s(0)>\max_i{f_i((\bbx(0),0))}$ so that $\bx(0) \in \widehat{\mathcal{D}}(0)$, i.e., the initial condition lies in the `enlarged' feasible set $\widehat{\mathcal{D}}(0)$. In the next lemma, we characterize the approximation error in terms of $c(t)$, $s(t)$, and the optimal dual variables in \eqref{eq: time_varying_KKT_conditions}.

\begin{lemma} \label{lemma: pertubation_suboptimality_bounds}
	Let $\bx^{\star}(t)$ and $\widehat{\bx}^{\star}(t)$ be defined as in \eqref{eq: time_varying_problem_with_barrier}  and \eqref{eq: time_varying_problem_with_log_barrier_and_slack}, respectively. Then, under Assumptions \ref{assumption: convexity} and \ref{assumption: strict_feasibility}, the following inequality holds for all $t\geq 0$,
	\begin{align} \label{eq: suboptimality_bound}
	|f_0(\widehat{\bbx}^{\star}(t),t)-f_0({\bx}^{\star}(t),t)| \leq \dfrac{p}{c(t)}+\sum_{i=1}^{p} \lambda^{\star}_i(t)s(t).
	\end{align}
\end{lemma}
\begin{myproof} 
See Appendix \ref{ap_perturbation_suboptimality_bound}.
\end{myproof}

The above lemma suggests that, if $s(t)$ and $c(t)$ are chosen such that the right-hand side of  $\eqref{eq: suboptimality_bound}$ converges to zero, the approximate solution $\widehat{\bx}^\star(t)$ converges to the optimal solution $\bx^\star(t)$ in \eqref{eq: time_varying_problem_with_barrier}.
In what follows, we design a dynamical system whose solution globally asymptotically converges to $\widehat{\bx}^\star(t)$.
Let us define ${\Phi}(\bx,c,s,t)$ as,
\begin{align} \label{eq: perturbed_time_varying_barrier_function}
\Phi(\bx,c,s,t) := f_0(\bx,t) - \frac{1}{c}\sum_{i=1}^p \log\left(s-f_i(\bx,t)\right).
\end{align}
Then, the optimal solution $\widehat{\bx}^\star(t)$ in \eqref{eq: time_varying_problem_with_log_barrier_and_slack} satisfies the optimality condition $\nabla_{\bx} \Phi(\widehat{\bx}^\star(t),c(t),s(t),t)=\mathbf{0}_n$ for all $t \geq 0$. Given the time evolution of $\bx(t)$, $s(t)$, and $c(t)$, we can use the chain rule to compute the time derivative of the gradient $\nabla_{\bx}{\Phi}(\bx,s,c,t)$, as follows:
\begin{align}
\dot\nabla_{\bx} {\Phi} = \nabla_{\bx\bx}\Phi \ \dot{\bx} + \nabla_{\bx s}\Phi \ \dot{s} + \nabla_{\bx c }\Phi \ \dot{c} + \nabla_{\bx t }\Phi.
\end{align}
Our goal is to design the dynamics of $\bx(t)$, as well as the time evolution of $c(t)$ and $s(t)$, such that the gradient $\nabla_{\bx} {\Phi}$ vanishes exponentially. In particular, based on the results in Section \ref{sec_main_section}, we propose the following dynamical system:
\begin{align} \label{eq: perturbed_time_varying_newton_barrier}
\dot{\bx}=-\nabla_{{\bx}{\bx}}^{-1}{\Phi}\Big[\, \alpha\nabla_{{\bx}}{\Phi}+\nabla_{\bx s}\Phi \ \dot{s} + \nabla_{\bx c }\Phi \ \dot{c} + \nabla_{\bx t }\Phi\, \Big].
\end{align}
The first term inside the bracket in the right-hand side corresponds to a Newton-like term, which is similar to the correction term in \eqref{eq: time_varying_newton}.
The remaining terms in \eqref{eq: perturbed_time_varying_newton_barrier} play a role similar to the prediction term in \eqref{eq: time_varying_newton}, since they account for time variations in $\Phi(\bx,c,s,t)$ through $s(t)$, $c(t)$, and $\{f_i(\bx,t)\}_{i=0}^{i=p}$.

Notice that it is important for the dynamics in \eqref{eq: perturbed_time_varying_newton_barrier} to render a solution $\bx(t)$ such that the argument of the logarithmic barrier functions in \eqref{eq: perturbed_time_varying_barrier_function} remains positive, i.e., we must have that $\bx(t) \in \widehat{\mathcal{D}}(t)$ for all $t\geq 0$. The next lemma states that this is, indeed, the case.
\begin{proposition} \label{lemma: time_varying_barrier_feasibility}
	Under Assumptions \ref{assumption: convexity}, \ref{assumption: uniform_strong_convexity}, and \ref{assumption: strict_feasibility}, the solution $\bx(t)$ to \eqref{eq: perturbed_time_varying_newton_barrier} satisfies $\bx(t) \in \widehat{\mathcal{D}}(t)$ for all $t\geq 0$ if $\bx(0) \in \widehat{\mathcal{D}}(0)$.
\end{proposition}
\begin{myproof}
	In Appendix \ref{ap_feasibility}.
\end{myproof}

In the following lemma, we prove that the solution to \eqref{eq: perturbed_time_varying_newton_barrier} converges exponentially to the approximate solution $\widehat{\bx}^{\star}(t)$ in \eqref{eq: time_varying_problem_with_log_barrier_and_slack}, for all initial conditions $\bx(0) \in \mathbb{R}^n$.
\begin{lemma}	\label{lemma: time_varying_barrier_convergence}
	Let $\widehat{\bx}^{\star}(t)$ be defined in \eqref{eq: time_varying_problem_with_log_barrier_and_slack} and ${\bx}(t)$ be the solution of \eqref{eq: perturbed_time_varying_newton_barrier} for  $\bx(0) \in \mathbb{R}^n$, $s(0)>\max_i{f_i((\bbx(0),0))}$, and $c(0)>0$. Then, under Assumptions \ref{assumption: convexity}, \ref{assumption: uniform_strong_convexity}, and \ref{assumption: strict_feasibility}, the following inequality holds,
	\begin{align} 
	\|{\bx}(t)-\widehat{\bx}^{\star}(t)\|_2\leq \dfrac{1}{m} \|\nabla_{\bx} {\Phi}(\bx(0),c(0),s(0),0)\|_{2} e^{-\alpha t}.\nonumber
	\end{align} 
\end{lemma}
	\begin{myproof}
		See Appendix \ref{ap: lemma_time_varying_barrier_convergence}. 
	\end{myproof}

Next, we need to establish the convergence of the approximate solution $\widehat{\bx}^\star(t)$ in \eqref{eq: time_varying_problem_with_log_barrier_and_slack} to the optimal solution $\bx^\star(t)$ in  \eqref{eq: time_varying_problem_with_barrier}. Intuitively, we need the barrier parameter $c(t)$ to asymptotically go to infinity and the slack variable $s(t)$ to asymptotically vanish so that the approximation error vanishes, according to \eqref{eq: suboptimality_bound}. For this to be true, we need to make the following assumption about the optimal dual variables defined in \eqref{eq: time_varying_KKT_conditions}.

\begin{assumption} \label{assumption: bounded_in_times}
	\normalfont For any $\gamma>0$, the optimal dual variables satisfy $\lambda_i^{\star}(t)\exp({-\gamma t})\to 0$ as $t \to \infty$, for all $i\in [p]$.
\end{assumption}
The above assumption excludes the possibility for the optimal dual variables to escape to infinity exponentially fast. By virtue of Assumption \ref{assumption: bounded_in_times}, the approximation error in  \eqref{eq: suboptimality_bound} vanishes asymptotically if the slack variable $s(t)$ goes to zero exponentially fast and the barrier parameter $c(t)$ diverges to infinity. The next theorem states the main result of this subsection.
\begin{theorem}\label{thm_main_theorem}
	Consider the optimization problem in \eqref{eq: time_varying_problem_with_barrier} and the objective function in \eqref{eq: perturbed_time_varying_barrier_function}. Let ${\bx}(t)$ be the solution of \eqref{eq: perturbed_time_varying_newton_barrier} with initial conditions $\bx(0) \in \mathbb{R}^n$, $c(0)>0$, and $s(0) > \max_{1 \leq i \leq p} f_i(\bx(0),0)$. Let $\lim\limits_{t \to \infty }c(t)=\infty$, and $s(t)=s(0)\exp({-\gamma_s t})$ for some $\gamma_s>0$. Then, under Assumptions \ref{assumption: convexity}, \ref{assumption: uniform_strong_convexity}, \ref{assumption: strict_feasibility}, and \ref{assumption: bounded_in_times}, we have that
	$$\lim\limits_{t \to \infty} \|{\bx}(t)-\bx^{\star}(t)\|_2 = 0.$$
\end{theorem}
\begin{myproof}
		By Assumption \ref{assumption: bounded_in_times} and Lemma \ref{lemma: pertubation_suboptimality_bounds}, we have that $\lim\limits_{t \to \infty} |f_0(\widehat{\bbx}^{\star}(t),t)-f_0({\bx}^{\star}(t),t)|=0$ when $s(t)=s(0)\exp(-\gamma_s t)$, and $\lim\limits_{t \to \infty} c(t)=\infty$. Strong convexity of $f_0$ (Assumption \ref{assumption: uniform_strong_convexity}) lets us to conclude that 	$\lim\limits_{t \to \infty} \|{\bx}(t)-\bx^{\star}(t)\|_2 = 0$.
\end{myproof}

According to Theorem \ref{thm_main_theorem}, the continuous-time dynamics in \eqref{eq: perturbed_time_varying_newton_barrier} yields a solution that asymptotically converges to the optimal solution in \eqref{eq: time_varying_problem_with_barrier} from arbitrary initial conditions. Some remarks are in order.

\begin{remark}[Barrier Parameter]
	{\rm The barrier parameter $c(t)$ is required to be positive, monotonically increasing, asymptotically converging to infinity, and bounded in finite time. A convenient choice is ${c}(t)=c(0) \exp({\gamma_c t})$ for $\gamma_c,c(0) > 0$. \ifx A more convenient choice is that $c(t)$ follows the adaptive rule \textcolor{red}{dot below?}
	\begin{align}
	\dot{c}(t)=\gamma_c . \Big[s(t)-\underset{1 \leq i \leq p}{\max} f_i(\bx(t),t)\Big].
	\end{align}
	\textcolor{red}{More clear... This particular choice does not increase the barrier parameter when the trajectory reaches the boundary. This has a practical importance because the gradient has rapid variations near the boundary and increasing $c$ when the trajectory is on the boundary might cause instability.}\fi}
\end{remark}

\begin{remark}[Adding Equality Constraints]
	{\rm As mentioned at the beginning of Section \ref{sec_interior_point}, we ignored equality constraints in our analysis. In order to account for equality constraints, we define the Lagrangian function as $\mathcal{L}(\bx,\bnu,c,s,t)=\Phi(\bx,c,s,t)+\bnu^\top(\mathbf{A}(t)\bx-\mathbf{b}(t))$ where $\Phi(\bx,c,s,t)$ is defined in \eqref{eq: perturbed_time_varying_barrier_function}, and consider the aggregate vector of decision variables $\bz=[\bx^\top \bnu^\top]^\top \in \mathbb{R}^{n+q}$. If we replace $\Phi$ by $\mathcal{L}$ in \eqref{eq: perturbed_time_varying_newton_barrier} and $\bx$ with $\bz$, we would obtain a dynamical system that solves the optimization problem in \eqref{eq: inequality_constrained_time_varying_problem}, where both equality and inequality constraints are considered. However, we need the inverse Hessian inverse in $\nabla_{\bz\bz}^{-1} \mathcal{L}$ to be uniformly bounded, i.e., there must exist an $M>0$ such that
\begin{align}
	\Big \|\begin{bmatrix}
	\nabla_{\bx\bx} \Phi(\bx,c,s,t) & \mathbf{A}^\top(t) \\ \mathbf{A}(t) & \mathbf{0}_{q \times q} 
	\end{bmatrix}^{-1}\Big \|_2  \leq M, \ \text{for all } t \geq 0.\nonumber
\end{align}
Therefore, using the same arguments as in Proposition \ref{lemma: equality_constrained_time_varyging_newton}, the result in Theorem \ref{thm_main_theorem} also holds, i.e., $\lim\limits_{t \to \infty} \|{\bz}(t)-\bz^{\star}(t)\|_2 = 0.$}
\end{remark}

%
\begin{remark}[Accelerating Standard Interior-Point Method]\label{remark: Acceleration of static interior point method}\normalfont As a particular application of our results, we consider the following time-invariant convex optimization problem,
	\begin{align} \label{eq: inequality_constrained_time_invariant_problem}
	\bbx^{\star}:=\arg\min  \ & f_{0}(\bbx) \ 
	\st \ f_{i}(\bbx) \leq 0, i\in [p].
	\end{align} 
	Using barrier functions, we define
	\begin{align} \label{eq: time_invariant_barrier}
	\Phi(\bbx,c)= f_0(\bbx) - \frac{1}{c}\sum_{i=1}^{p}\log(-f_i(\bbx)),
	\end{align}
	and the corresponding central path
	\begin{equation}\label{eqn_relaxed_time_invarant}
	\bbx^{\star}(c):=\arg\min_{\bbx\in \mathcal{D}} \, \Phi(\bbx,c),
	\end{equation}
	where $\mathcal{D}=\{\bbx\in \mbR^n \colon f_i(\bbx)<0, i\in [p]\}$ is the \emph{interior} of the feasible set. It follows from Lemma \ref{lemma: pertubation_suboptimality_bounds} that $\lim\limits_{c \to \infty}\|\bbx^{\star}(c)- \bbx^{\star}\|=0$. In the standard implementation of the interior-point method \cite[Chapter 11]{boyd2004convex}, the optimization problem \eqref{eqn_relaxed_time_invarant} is solved sequentially for a positive growing sequence $(c_k)_{k=1}^{\infty}$, each starting from the optimal
	solution of the previous optimization problem. The resulting sequence $(\bx^\star(c_k))_{k=1}^{\infty}$ converges to the optimal point $\bbx^{\star}$ as $c_k \to \infty$. For each fixed $c_k$, $\bx^{\star}(c_k)$ can be found, for instance, using the following Newton-like dynamics:
	\begin{align} \label{eq: time_invariant_newton_barrier}
	\dot \bx_k(t) = -\nabla_{\bx\bx}^{-1} \Phi(\bx_k(t),c_{k})\nabla_{\bx}\Phi(\bx_k(t),c_{k}).
	\end{align}
	According to Theorem \ref{thm_main_theorem},  we have that $\lim\limits_{t \to \infty} \|\bbx_{k}(t)-\bbx^\star(c_k)\|=0$. Thus, by choosing the initial points as $\bbx_{k}(0)=\lim\limits_{t \to \infty}\bbx_{k-1}(t)$, we can build a continuous path that converges to $\bbx^\star$ in \eqref{eq: inequality_constrained_time_invariant_problem}. As a less computationally expensive alternative, we propose to consider an increasing time-dependent barrier parameter $c=c(t)$, in lieu of discontinuous jumps. In this case, the problem in \eqref{eqn_relaxed_time_invarant} renders a time-varying optimization. Hence, the continuous dynamics \eqref{eq: perturbed_time_varying_newton_barrier} using the barrier function in \eqref{eq: time_invariant_barrier} yields a solution $\bx(t)$ satisfying $\lim\limits_{t \to \infty}\|\bx(t)-\bx^\star\|_2=0$ when $\lim\limits_{t \to \infty} c(t) = \infty$. We will numerically illustrate the performance of this approach in Subsection \ref{sec_static_interior_point}.
\end{remark}

%
%
\begin{remark}[Second-Order Implementation]\normalfont
Logarithmic barrier functions present a singularity at the origin that may induce numerical instability in the discrete-time implementation of the dynamical system \eqref{eq: perturbed_time_varying_newton_barrier}. To avoid this issue, we replace the first-order dynamics in \eqref{eq: perturbed_time_varying_newton_barrier} by the following second-order dynamics \cite{su2009traffic}:
	\begin{align} \label{eq: second-order dynamics}
	\dot{\bx}&=-\nabla_{{\bx}{\bx}}^{-1}{\Phi} \Big[\, \alpha \by+ \nabla_{\bx s}\Phi \ \dot{s} + \nabla_{\bx c }\Phi \ \dot{c}+\nabla_{{\bx} t}{\Phi}\, \Big], \nonumber \\
	\dot{\by}&=-\gamma \by + \alpha\nabla_{{\bx}}{\Phi},
	\end{align}
	where $\gamma >0$ is an arbitrary constant. Intuitively, the gradient function is passed through a first-order low-pass filter (the second ODE in \eqref{eq: second-order dynamics}) whose output is then fed into the main dynamics. The resulting dynamics tends to reduce numerical instability induced by discretization. It was shown in \cite{su2009traffic} that $V=\frac{1}{2} \nabla_{\bx} \Phi^\top \nabla_{\bx} \Phi+ \frac{1}{2} \by^\top \by$ is a Lyapunov function for \eqref{eq: second-order dynamics}, proving that $\lim\limits_{t \to \infty} \|\nabla_{\bx} \Phi\|=0$. 
\end{remark}

%
\subsection{Online Implementation}
The dynamical system proposed in \eqref{eq: perturbed_time_varying_newton_barrier} includes the prediction term $\nabla_{\bx t} \Phi(\bx,t)$, whose computation involves finding the terms $\{\nabla_{\bx t} f_i(\bx,t)\}_{i=0}^{p}$ and $\{ \frac{\partial}{\partial t}f_i(\bx,t)\}_{i=1}^{p}$. In an online setting, we might only have access to limited or noisy information about these terms. More precisely, assume that we have access to an estimate of $\nabla_{\bx t} \Phi$ denoted by $\widehat{\nabla}_{\bx t} \Phi$ that satisfies the bound
\begin{align} \label{eq: prediction_dynamics_estimate}
\|\widehat{\nabla}_{\bbx t} \Phi- {\nabla}_{\bbx t} \Phi\|_2 \leq \eta,
\end{align}
for some known $\eta>0$. In this setting, we consider the following dynamics:
\begin{align} \label{eq: perturbed_time_varying_newton_barrier_robust}
\dot{\bx}=-\nabla_{{\bx}{\bx}}^{-1}{\Phi}\Big[\, \alpha \nabla_{{\bx}}{\Phi}+\nabla_{\bx s}\Phi \ \dot{s} + \nabla_{\bx c }\Phi \ \dot{c} + \widehat{\nabla}_{\bx t }\Phi\, \Big],
\end{align}
where we define a state-dependent $\alpha=\alpha(\bx)$ as follows:
\begin{align} \label{eq: adaptive alpha}
\alpha=\frac{\alpha_0}{\max(\|\nabla_{\bbx} \Phi\|_2,\varepsilon)},
\end{align}
where $\alpha_0$ satisfies $\alpha_0>\eta$, and $\varepsilon>0$ is an arbitrary constant. The next theorem states that the solution of \eqref{eq: perturbed_time_varying_newton_barrier_robust} converges to an $\varepsilon$-neighborhood of the approximate optimal solution $\widehat{\bx}^\star(t)$, defined in \eqref{eq: time_varying_problem_with_log_barrier_and_slack}, in finite time and will stay there forever.
\begin{theorem} \label{lem: unconstrained_finite_time}
	Denote $\bx(t)$ as the solution of \eqref{eq: perturbed_time_varying_newton_barrier_robust} where $\Phi$ is defined in \eqref{eq: perturbed_time_varying_barrier_function}. Assume $\widehat{\nabla}_{\bx t} \Phi$ satisfies the bound in \eqref{eq: prediction_dynamics_estimate}, and the coefficient $\alpha$ is defined in \eqref{eq: adaptive alpha} with $\alpha_0>\eta$ and $\varepsilon>0$. Then, under Assumptions \ref{assumption: convexity}, \ref{assumption: uniform_strong_convexity}, and \ref{assumption: strict_feasibility}, the solution $\bx(t)$ converges to the set $\mathcal{S}_t({\varepsilon}):=\{\bx \in \mathbb{R}^n \colon \|\nabla_{\bx} \Phi(\bx,s(t),c(t),t)\|_2 \leq \varepsilon\}$ in finite time.
\end{theorem}

\begin{myproof}
See Appendix \ref{ap: lem: unconstrained_finite_time}.
\end{myproof}
%
%


\section{Numerical Experiments}
In this section, we provide three numerical examples to illustrate the time-varying optimization framework herein proposed. In Subsection \ref{sec_sim_sinthetic}, we solve a synthetic optimization problem to illustrate the effectiveness of the prediction-correction interior-point method in solving inequality-constrained problems. In Subsection \ref{sec_static_interior_point}, we use the accelerated interior-point method discussed in Remark \ref{remark: Acceleration of static interior point method} to solve a large-scale $\ell_1$-regularized least-squares problem. In Subsection \ref{sec_robot_navigation}, we solve a navigation problem to drive a disk-shaped robot towards a (potentially moving) desired location $\bbx_d$ without colliding with obstacles in the environment.

\subsection{Time-Varying Quadratic Programming}\label{sec_sim_sinthetic} 
Consider the following TV quadratic optimization problem:
    \begin{alignat}{2} \label{eq: example_quadtaric}
    \bbx^{\star}(t)
    :=\  &\argmin\ && \frac{1}{2} (x_1+\sin(t))^2+\frac{3}{2} (x_2+\cos(t))^2,  \nonumber \\ 
    &\st    \ && x_2-x_1 -\cos(t) \leq 0.
    \end{alignat} 
In the following simulation, we show how to track $\bx^\star(t)$ using the continuous-time dynamics in \eqref{eq: perturbed_time_varying_newton_barrier}. In order to illustrate the usage of the time-dependent slack variable $s(t)$, we choose the initial condition $\bx(0)=(-2,0)^\top$, which is infeasible at $t=0$. As discussed in Subsection \ref{sec_interior_point}, we include the slack variable to enlarge the feasible set. In this example, the augmented objective function in \eqref{eq: perturbed_time_varying_barrier_function} takes the form:
\begin{align*}
\Phi(\bx,s,c,t) =& \frac{1}{2} \Big(x_1+\sin(t)\Big)^2+\frac{3}{2} \Big(x_2+\cos(t)\Big)^2 \nonumber \\&- \frac{1}{c} \log\Big(s+\cos(t)+x_1-x_2\Big).
\end{align*}
\begin{figure}
    \centering
    \includegraphics[width=0.45\textwidth]{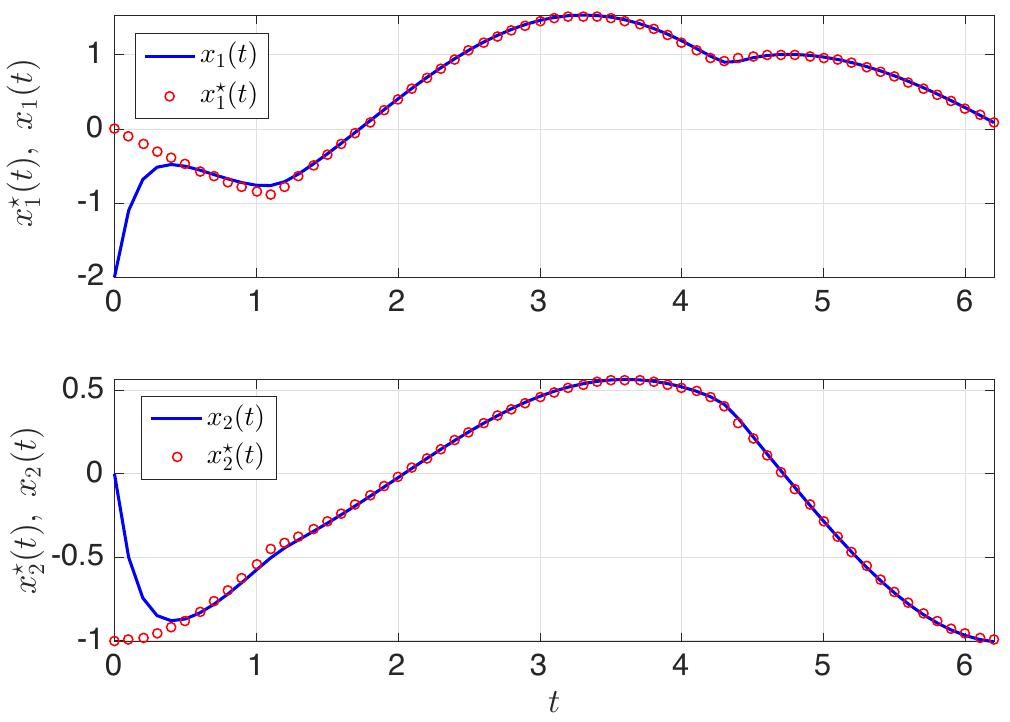}
    \caption{{\small Plot of the coordinates of the optimal trajectory $\bx^\star(t)$, defined in \eqref{eq: example_quadtaric}, and the tracking trajectory $\bx(t)$.}}
    \label{fig: example_quadtaric_trajectory}
\end{figure}

\begin{figure}
    \centering
    \includegraphics[width=0.47\textwidth]{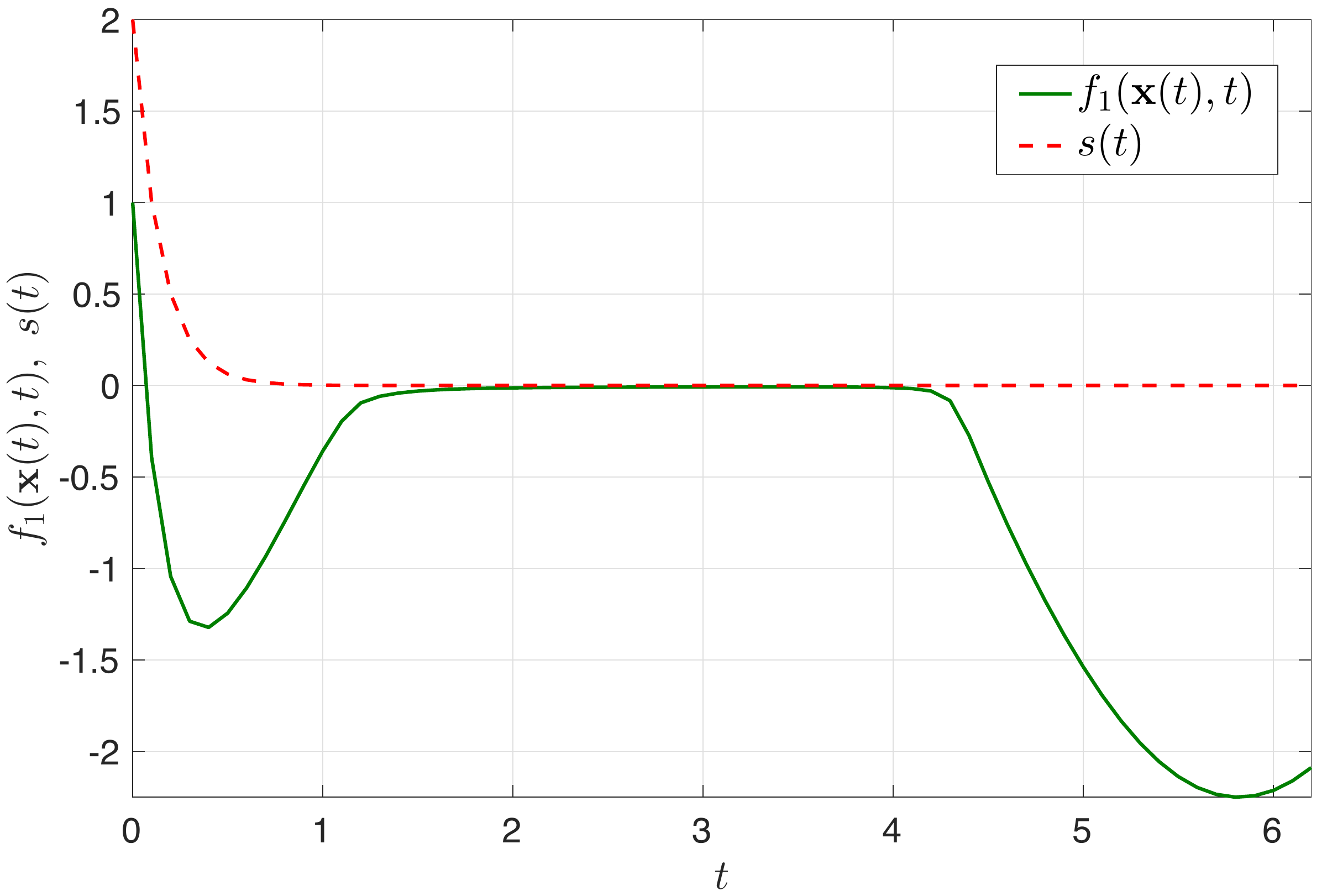}
    \caption{{\small Plot of the constraint function $f_1(\bx(t),t)=x_2(t)-x_1(t)-\sin(t)$ and the slack variable $s(t)$ against $t$.}}
    \label{fig: example_quadtaric_constraint_violation}
\end{figure}

In our simulations, we consider the following time-dependent barrier parameter and slack variable: $c(t)=10e^t$ and $s(t)=2e^{-5t}$. The initial value of the slack variable is $s(0)=2$; hence, $\bx(0)$ is initially feasible with respect to the enlarged feasible set. Using these particular values, all the conditions of Theorem~\ref{thm_main_theorem} are satisfied. We numerically solve the ODE in \eqref{eq: perturbed_time_varying_newton_barrier} for the time interval $t\in [0,2\pi]$ using Euler's forward method with step size $\tau=0.1$. In Fig. \ref{fig: example_quadtaric_trajectory}, we plot the trajectory of the resulting solution $\bx(t)=(x_1(t),x_2(t))^\top$ along with the optimal solution $\bx^\star(t)=(x^\star_1(t),x^\star_2(t))$ defined in \eqref{eq: example_quadtaric}. In Fig. \ref{fig: example_quadtaric_constraint_violation} we plot the time evolution of the constraint function $f_1(\bx,t):=x_2(t)-x_1(t) -\cos(t)$, as well as the slack variable $s(t)$.  Notice how, at $t=0$, the state $\bx(t)$ violates the constraint $f_1(\bx,t)\leq 0$. However, $\bx(t)$ converges to the feasible set exponentially fast as the slack variable $s(t)$ vanishes exponentially. \ifx In Fig. \ref{fig: example_quadtaric_alpha}, we plot the time evolution of the tracking error $\|\bx(t)-\bx^\star(t)\|_2$ for various values of $\alpha$ DEFINE WHERE? WHAT IS THE MEANING?. Observe how, as the value of $\alpha$ increases, the solution converges faster to the optimal trajectory and the amplitude of spikes in the tracking error decreases.\fi

\ifx
\begin{figure}
    \centering
    \includegraphics[width=0.5\textwidth]{Example_1_Fig_4}
    \caption{{\small Plot of the tracking error $\|\bx(t)-\bx^\star(t)\|_2$ as a function of $t$, and for various values of $\alpha$. Higher values of $\alpha$ results in faster convergence toward the optimal trajectory.}}
    \label{fig: example_quadtaric_alpha}
\end{figure}
\fi

\subsection{$\ell_1$-Regularized Least Squares}\label{sec_static_interior_point}
In this subsection, we illustrate how to use the accelerated interior-point method (described in Remark \ref{remark: Acceleration of static interior point method}) to solve the following $\ell_1$-regularized least-squares problem:
\begin{align} \label{eq: L1_LS}
\bx^\star (\lambda) = \underset{\bx \in \mathbb{R}^n}{\argmin} \ \|\mathbf{A}\bx-\mathbf{b}\|_2^2 + \lambda \|\bx\|_1,
\end{align}
where $\mathbf{A} \in \mathbb{R}^{m\times n} $ and $\mathbf{b} \in \mathbb{R}^m$ are given, and  $\lambda >0$ is a regularizer which is commonly used to prevent overfitting whenever $m<n$. Since the objective function in \eqref{eq: L1_LS} is not differentiable, we analyze the following (differentiable) equivalent convex program \cite{kim2007interior}:
\begin{align} \label{eq: L1_LS_differentiable}
(\bx^\star (\lambda),\bu^\star(\lambda) )=\  &\underset{\bx,\bu  \in \mathbb{R}^n}{\argmin} \ &&\|\mathbf{A}\bx-\mathbf{b}\|_2^2 + \lambda \sum_{i=1}^{n} u_i,  \nonumber \\ 
&\st    \ && -u_i \leq x_i \leq u_i, \ i \in [n].
\end{align}
%

\begin{figure}
    \centering
    \includegraphics[width=0.5\textwidth]{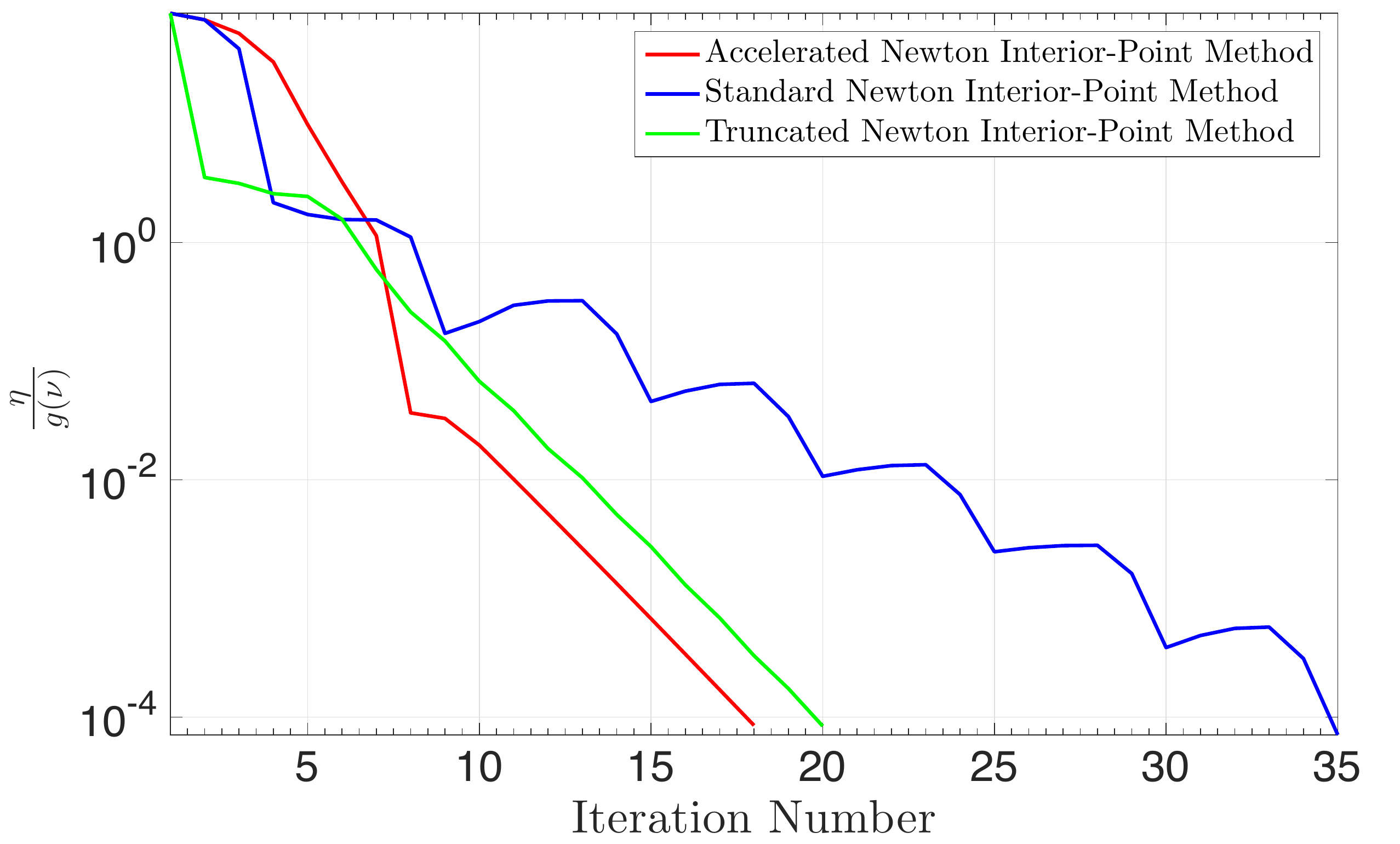}
    \caption{{\small Plot of the upper bound $\eta/g(\nu)$ on the relative suboptimality against the number of iterations for the three algorithms described in Subsection \ref{sec_static_interior_point}.}}
    \label{fig: example_2_LS}
\end{figure}

In the following numerical experiment, we generate a sparse vector $\bx^\star \in \mathbb{R}^{2048}$ with 10 entries of value $\pm 1$, and all other entries equal to zero. The entries of the measurement matrix $\mathbf{A} \in \mathbb{R}^{256\times 1024}$ are independently generated according to the standard normal density. The measured vector $\mathbf{b} \in \mathbb{R}^{256}$ is generated by $\mathbf{b}=\mathbf{A}\bx^\star+\bv$ where $\bv$ is a contaminating noise drawn from the Gaussian distribution $\mathcal{N}(\mathbf{0}_{256}, 0.01\, \mathbf{I}_{256})$. The regularizer parameter is chosen to be $\lambda=2$. For these numerical values, we solve \eqref{eq: L1_LS_differentiable} using three methods: (\emph{i}) the Standard Newton Interior-Point Method (SNIPM) \cite[Chapter 11]{boyd2004convex} where the central points $\bx^\star(c)$ are computed using Newton's method with the sequence $c_k=10 \times 5^k,\ k=0,1,\ldots$ ; (\emph{ii}) the Accelerated Newton Interior-Point Method (ANIPM), described in Remark \ref{remark: Acceleration of static interior point method}, where the barrier parameter is equal to $c(t)=10e^t$;
and (\emph{iii}) the Truncated Newton Interior-Point Method (TNIPM), described in \cite{kim2007interior}, where a preconditioned conjugate gradient method was proposed to compute the Newton step, and the barrier parameter is updated at each iteration. For all these three methods, we use a backtracking line search to adaptively select the step size. To assess the progress of the algorithms, we use the following quantity (as proposed in \cite{kim2007interior}),
\begin{align} \label{eq: LS_Example_relative_gap}
\frac{\eta}{g(\bnu)} := \frac{\|\mathbf{A}\bx-\mathbf{b}\|_2^2 + \lambda \|\bx\|_1 - g(\bnu)}{g(\bnu)}.
\end{align}
Here, $g(\bnu)$ is the dual function of the constrained problem akin to \eqref{eq: L1_LS}:
\begin{align}
\min \ \bz^\top \bz+ \lambda \|\bx\|_1 \mbox{  such that } \ \bz=\mathbf{A} \bx - \mathbf{b},
\end{align}
and $\bnu \in \mathbb{R}^m$ is the dual vector associated with the constraint $\bz=\mathbf{A} \bx - \mathbf{b}$. The quantity in \eqref{eq: LS_Example_relative_gap} is an upper bound of the relative duality gap $[p^\star-g(\bnu)]/g(\bnu)$, where $p^\star=\|\mathbf{A}\bx^\star(\lambda)-\mathbf{b}\|_2^2 + \lambda \|\bx^\star(\lambda)\|_1$ is the primal optimal value (see \cite{kim2007interior} for more details). Fig. \ref{fig: example_2_LS} illustrates the evolution of $\eta/g(\bnu)$ against the iteration number for these three algorithms. For the stopping criterion, we choose $\eta/g(\bnu) \leq 10^{-4}$. In our simulations, the SNIPM takes $42$ iterations, while the ANIPM proposed in Remark \ref{remark: Acceleration of static interior point method} takes $17$ iterations. Notice that the performance of our accelerated method is comparable to TNIPM, since in the latter method the barrier parameter is also updated at each iteration, but the prediction term is not included.

\subsection{Robot Navigation}\label{sec_robot_navigation}
In this subsection, we solve the navigation problem of driving a disk-shaped robot of radius $r>0$ to a given configuration $\bbx_d$ without colliding with obstacles in the environment. More precisely, let us consider a closed and convex workspace $\ccalW \subset \mathbb{R}^n$ of possible configurations that the robot can take. Assume that the workspace is populated with $m$ non-intersecting spherical obstacles, where the center and radius of the $i$-th obstacle are denoted by $\bbx_i \in \ccalW$ and $r_i >0$, respectively. We define the \emph{free space}, denoted by \ccalF, as the set of configurations in the workspace in which the robot does not collide with any of the obstacles. More formally,
\begin{equation}
\ccalF = \left\{\bbx \in \ccalW \colon \overline{B}(\bbx,r)\subseteq \mathcal{W} \setminus \cup_{i=1}^m B(\bbx_i,r_i ) \right\},
\end{equation} 
where $B(\bbx,r)$ is the $n$-dimensional open ball centered at $\bbx$ with radius $r$, and $\overline{B}(\bbx,r)$ represents its closure.

Let us denote the center of mass of the robot by $\bbx_c$.  Given a final desired configuration $\bbx_d\in\ccalF$, the navigation problem under consideration consists of finding a trajectory of $\bbx_c$ such that $\bbx_c(t) \in \ccalF$ for all $t\geq 0$, and $\lim_{t\to\infty} \bbx_c(t) = \bbx_d$.
In \cite{arslan_kod_ICRA2016B}, the authors proposed a solution to this problem using the idea of \emph{projected goal}, as described below. This idea consists of continuously computing the projection of the destination $\bbx_d$ onto a neighborhood around the center of mass of the robot in which there are no obstacles. Denote this projection by $\bar{\bbx}_d$ ---yet to be properly defined--- then, the control law $\dot{\bbx}_c=-\left(\bbx_c-\bar{\bbx}_d\right)$ ensures convergence of the center of mass of the robot to the desired configuration while avoiding the obstacles \cite{arslan_kod_ICRA2016B}. As we describe below, this technique can be interpreted as the solution of a TV convex optimization problem. To formulate this problem, we first need to provide some definitions.      

\begin{figure}
\centering
\includegraphics[width=0.30\textwidth]{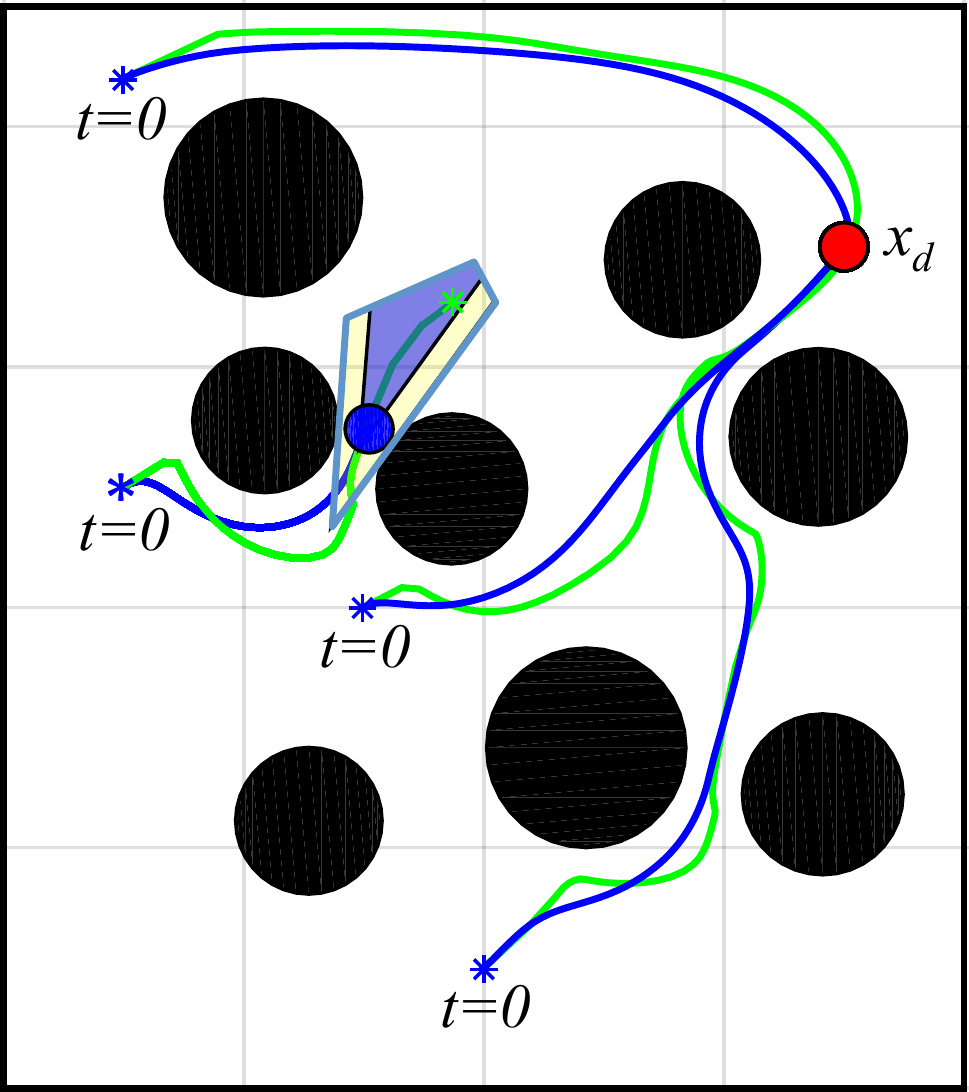}
\caption{{\small The red circle represents a desired configuration $\bbx_d$. The green and blue lines represent, respectively, the trajectories of the estimates of the projected goal $\hat{\bbx}(t)$ and the trajectories of the robot $\bbx_c(t)$ for 4 different initial conditions.}}
\label{fig_static_goal_trajectory}
\end{figure}

We define first the notion of \emph{power distance} between a point $\bbx$ and a disk $B(\bbx_i,r_i )$ as $\mathcal{P}(\bbx,B(\bbx_i,r_i ))=\left\Vert\bbx-\bbx_i\right\Vert_2^2-r_i^2$, \cite{aurenhammer1987power}. We define the so-called \emph{local workspace} around $\bbx_c$ as
\[
\ccalL \ccalW(\bbx_c)=\left\{\bbx \in \ccalW \colon \mathcal{P}(\bbx,B(\bbx_c,r))\leq \mathcal{P}(\bbx,B(\bbx_i,r_i )),\, \forall i \right\},
\]
i.e., the set of points in $\mathcal{W}$ that are closer (in power distance) to the robot than to any of the obstacles.
The local workspace defines a polytope whose boundaries are hyperplanes, such as the polygon marked with a thick light blue line in Fig. \ref{fig_static_goal_trajectory} (see \cite[Eq. (6)]{arslan_kod_ICRA2016B} for an explicit expression of these hyperplanes).
Furthermore, the \emph{collision-free local workspace} around $\bbx_c$ is defined as \cite{arslan_kod_ICRA2016B}:
\begin{equation}
\ccalL \ccalF (\bbx_c) = \left\{ \bbx\in \ccalW \colon \bba_i(\bbx_c)^\top \bbx - b_i(\bbx_c) \leq 0, \; i=1\ldots m\right\},\nonumber
\end{equation}
where,
\begin{align}\label{eqn_bi}
\bba_i(\bbx_c) & = \bbx_i - \bbx_c, \quad \theta_i(\bbx_c) = \frac{1}{2}     - \frac{r_i^2 - r^2}{2\| \bbx_i - \bbx_c\|^2}, \\
b_i(\bbx_c) & = (\bbx_i - \bbx_c)^\top\left( \theta_i\bbx_i +(1-\theta_i)\bbx_c+r\frac{\bbx_c-\bbx_i}{\|\bbx_c-\bbx_i\|} \right) \nonumber.
\end{align}
Assuming that the robot follows the integrator dynamics $\dot{\bbx}_c= \bbu(\bbx_c)$, the controller proposed in \cite{arslan_kod_ICRA2016B} is given by 
\begin{equation}\label{eqn_omur_controller}
\dot{\bbx}_c = -K ( \bbx_c - \bbx^\star),
\end{equation}
where $K>0$ is the gain of the controller and $\bbx^\star$ is the orthogonal projection of the desired configuration $\bbx_d$ onto the collision-free local workspace $\ccalL \ccalF(\bbx_c)$. Under the assumption that the distance between the center of any two obstacles $i$ and $j$ is larger than $r_i+r_j +2r$, it can be shown that the controller law in \eqref{eqn_omur_controller} solves the navigation problem (\cite[~Theorem 1]{arslan_kod_ICRA2016B}). In what follows, we cast the navigation problem as a TV convex optimization program that can be solved using the tools developed in this paper.

\begin{figure*}
     \centering
    \includegraphics[width=1.0\textwidth]{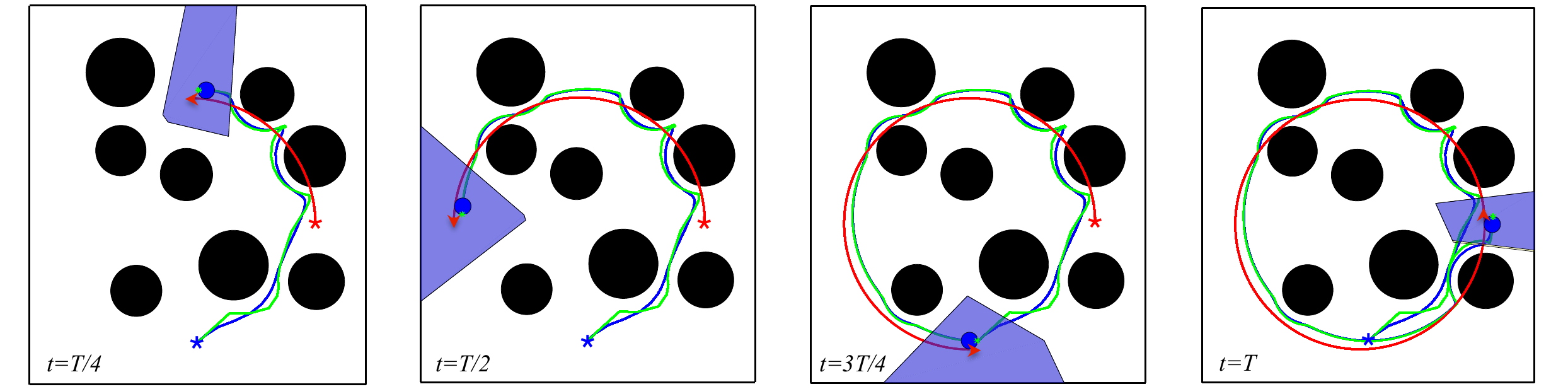}
        ~ \caption{{\small  Trajectory of the  estimate of the projected goal $\hat{\bbx}$ (the green line, starting at the blue star) and  trajectory of the robot (the blue line, starting at the blue star) tracking a moving target following a circular trajectory (red line, starting at the red star). The subplots correspond to time instances $t=T/4,T/2,3T/4,$ and $T$. The local workspace for each time instance is represented by a solid purple polygon containing the disk-shaped robot (the blue circle).}}\label{fig_moving}
\end{figure*}

\subsubsection{Interior-Point Method to Estimate the Projected Goal}\label{sec_static_target}
We now show that the prediction-correction interior-point method developed in Section \ref{sec_main_section} can be used to efficiently build an estimator $\hat{\bbx}$ of the projection of $\bbx_d$ onto the set $\ccalL \ccalF(\bbx_c)$, which we denote by $\bbx^\star$. First, observe that $\bbx^\star$ can be defined as the solution of the following convex optimization problem,
\begin{align}\label{eqn_projected_goal}
\bbx^\star:= &\argmin_{\bbx \in \mathbb{R}^n} \frac{1}{2} \| \bbx-\bbx_d\|^2 \nonumber \\
&\mbox{s.t.} \quad \bba_i(\bbx_c)^\top\bbx - b_i(\bbx_c) \leq 0, \quad i =1\ldots m.
\end{align}
Observe that since $\bba_i(\bbx_c)$ and $b_i(\bbx_c)$ depend on the position of the center of mass of the robot, the above optimization problem has an implicit dependence on time through $\bx_c$. We estimate the projected goal $\bx^\star$ as the solution to the ODE in \eqref{eq: perturbed_time_varying_newton_barrier} with initial condition ${\bbx}(0) = \bbx_c(0)$, i.e., the initial position of the robot, and the following objective function:
\begin{equation}
\Phi(\bbx,\bbx_c,t) = \frac{1}{2} \|\bbx-\bbx_d\|^2 - \frac{1}{c(t)}\sum_{i=1}^m \log(b_i(\bbx_c) - \bba_i(\bbx_c)^\top\bbx).\nonumber
\end{equation}
In Appendix \ref{ap_numerical_examples}, we derive explicit expressions for all the terms involved in this ODE.

Next, we consider the control law \eqref{eqn_omur_controller}, with the difference that we use an estimate of the projected goal instead of the projected goal itself, i.e., we consider the closed loop dynamics
\begin{equation}\label{eqn_control_law}
\dot{\bbx}_c = -K(\bbx_c-\hat{\bbx}),
\end{equation}
where the estimator $\hat{\bbx}(t)$ is the solution to the ODE in \eqref{eq: perturbed_time_varying_newton_barrier} with initial condition $\hat{\bbx}(0) = \bbx_c(0)$. An important feature of $\hat{\bbx}(t)$ is that it is feasible at all times, i.e., $\hat{\bbx}(t) \in \ccalL \ccalF(\bbx_c)$. This follows from Lemma~\ref{lemma: time_varying_barrier_feasibility} and the fact that the initial condition is assumed to be feasible, i.e., the robot is initially located in the free space. Moreover, the estimator $\hat{\bbx}(t)$ converges exponentially to the projection of $\bbx_d$ onto the collision-free local workspace, denoted by $\bbx^\star$.

%
%

To evaluate the performance of the proposed controller and optimizer, we consider a square workspace $\ccalW=[-20,20]^2$ containing  $8$ circular obstacles (black circles in Fig. \ref{fig_static_goal_trajectory}). In Fig. \ref{fig_static_goal_trajectory}, we also depict the trajectories followed by a disc-shaped robot of radius equal to one (blue circle) for four different initial conditions.
The green and blue lines represent, respectively, the trajectories of the estimates  $\hat{\bbx}(t)$ of the projected goal onto the collision-free local workspace, and the trajectories of the center of mass of the robot $\bbx_c(t)$ for 4 different initial conditions. The blue circle represents a particular configuration of the robot, where the local workspace $\ccalL \ccalW (\bbx)$ (resp., the collision-free local workspace $\ccalL \ccalF (\bbx_c)$) is the polygon enclosed within light blue lines (resp., the polygon filled in solid purple).
For these particular realizations, we have set $\alpha=5$ in \eqref{eq: perturbed_time_varying_newton_barrier}, and $K=0.01$ in \eqref{eqn_control_law}. Finally, the barrier parameter in \eqref{eq: perturbed_time_varying_barrier_function} is chosen to be $c(t) = e^{0.001 t}$. In Fig. \ref{fig_static_goal_trajectory}, we observe how the robot succeeds in converging to the desired destination. Collision avoidance is ensured due to the fact that the estimate of the projected goal $\hat{\bbx}$ remains always in the collision-free local workspace $\ccalL\ccalF(\bbx_c)$.
%
%

\subsubsection{Moving Targets}\label{sec_moving_target}
In our final experiment, we consider a similar navigation problem, but in this case the disk-shaped robot must track a moving target, i.e., $\bbx_d(t) : \mathbb{R}_+ \to \ccalW$. In this scenario, no theoretical guarantees are provided by the technique proposed in \cite{arslan_kod_ICRA2016B}; however, we demonstrate how our approach can be used to track a moving target.

In our experiment, we consider a moving target $\bbx_d(t)$ following a circumference of radius $15$, centered at the origin, and moving periodically with a period $T= 2\times10^{3} s$. Observe that the target trajectory (red line in Fig. \ref{fig_moving}) is allowed to intersect the circular obstacles (e.g., obstacles are on the ground, and the target is flying).
To track this target, we use the controller \eqref{eqn_control_law} where $\hat{\bbx}$ is the solution to the ODE in \eqref{eq: perturbed_time_varying_newton_barrier} with the following barrier function:
\begin{equation}
\Phi(\bbx,\bbx_c,t) = \frac{1}{2} \|\bbx-\bbx_d(t)\|^2 - \frac{1}{c(t)}\sum_{i=1}^m \log(b_i(\bbx_c) - \bba_i(\bbx_c)^\top\bbx),\nonumber
\end{equation}
where $\bba_i(\bbx_c) = \bbx_i-\bbx_c$ and $b_i(\bbx_c)$ are given by \eqref{eqn_bi}.  The parameter selection for our simulation is $K=0.05$, $\alpha = 30$, and $c(t) = 100e^{\alpha t}$ with $\alpha =0.001$. In Fig. \ref{fig_moving}, we depict the trajectory followed by the disk-shaped robot (blue circle) over time.  As we can observe, the robot succeeds in tracking the moving target while avoiding the circular obstacles.



\section{Conclusions}\label{se:conclusions}
In this paper, we have developed a prediction-correction interior-point method for solving convex optimization problems with time-varying objective and/or constraint functions. Using log-barrier penalty functions, we have proposed a continuous-time dynamical system for tracking the solution of the time-varying problem. This dynamical system contains both a correction term, which is a continuous-time implementation of Newton's method, as well as a prediction term that takes into account the time-varying nature of the objective and constraint functions. Under reasonable assumptions, our method globally asymptotically converges to the time-varying optimal solution of the original problem with a vanishing tracking error.
We have illustrated the applicability of the proposed method in two practical applications: a sparsity promoting least squares problem and a collision-free robot navigation problem.

\appendix
\section{Appendix}

	\subsection{Proof of Proposition \ref{lem: Time_Varying_Newton}} \label{ap: time_varying_unconstrained_newton}
	
		First, notice that by Assumption \ref{assumption: uniform_strong_convexity} (uniform strong convexity), the Hessian inverse ${\nabla_{\bx\bx}^{-1} f_0(\bx(t),t)}$ is defined and uniformly bounded for all $t \geq 0$. The time derivative of the gradient at $(\bx(t),t)$ can be written as 
		\begin{align} \label{eq: gradient_dynamics}
		\dot{\nabla}_{\bx} f_0(\bx(t),t)= \nabla_{\bx\bx} f_0(\bx(t),t) \dot{\bx}(t)+ \nabla_{\bx t} f_0(\bx(t),t). 
		\end{align}
		Substituting $\dot{\bx}(t)$ from \eqref{eq: time_varying_newton} in \eqref{eq: gradient_dynamics}, we obtain
		\begin{align*}
		 \dot\nabla_{\bx} f_0(\bx(t),t)=-\alpha \nabla_{\bx} f_0(\bx(t),t). 
		\end{align*}
		This is a first-order linear differential equation on $\nabla_{\bx}f_0(\bx(t),t)$, yielding the solution
		\begin{align*} \label{eq: unconstrained_time_varying_}
		\nabla_{\bx} f_0(\bx(t),t)=e^{-\alpha t}\nabla_{\bx} f_0(\bx(0),0),
		\end{align*}
		where $\bx(0) \in \mathbb{R}^n$ is the initial point. Apply Cauchy-Schwartz inequality, we obtain
		\begin{equation}\label{eqn_lemma1_aux}
		\|\nabla_{\bx} f_0(\bx(t),t)\|_2  \leq e^{-\alpha t} \|\nabla_{\bx} f_0(\bx(0),0)\|_2.
		\end{equation}
		Next, we fix a $t \geq 0$ and use the mean-value theorem to expand $\nabla_{\bx} f_0(\bx(t),t)$ around the optimal point $\nabla_{\bx} f_0(\bx^\star(t),t)=\mathbf{0}$,
		\begin{align}
		\nabla_{\bx} f_0(\bx(t),t) =\nabla_{\bx\bx} f_0(\boldsymbol{\eta}(t),t) (\bx(t)-\bx^\star(t)),\nonumber
		\end{align}
		where $\boldsymbol{\eta}(t)$ is a convex combination of $\bx(t)$ and $\bx^\star(t)$. It follows from uniform strong convexity of $f_0(\bx,t)$ (Assumption \ref{assumption: uniform_strong_convexity}) that $\|\nabla_{\bx\bx}^{-1} f_0(\bz(t),t) \|_2 \leq m^{-1}$. Whence, we can invoke \eqref{eqn_lemma1_aux} to write
		\begin{align*} 
		\|\bx(t)-\bx^{\star}(t)\|_2  & =\|\nabla_{\bx\bx}^{-1} f_0(\boldeta(t),t)\nabla_{\bx}f_0(\bx(t),t)\|_2 \\
		&\leq \|\nabla_{\bx\bx}^{-1} f_0(\boldeta(t),t)\|_2 \|\nabla_{\bx}f_0(\bx(t),t)\|_2 \\
		& = \dfrac{1}{m}\|\nabla_{\bx} f_0(\bx(0),0)\|_2 e^{-\alpha t}.
		\end{align*}
		The proof is complete.
		\subsection{Proof of Proposition \ref{lemma: equality_constrained_time_varyging_newton}} \label{ap: time_varying_equality_newton}

		The Hessian of the Lagrangian in \eqref{eq: time_varying_equality_lagrangian} with respect to $\bz=[\bx^\top, \ \bnu^\top]^\top$ is given by
		\begin{align*}
		\nabla_{\bz\bz} {\mathcal{L}}(\bz,t)=\begin{bmatrix} \nabla_{\bx\bx} f_0(\bx,t) & \mathbf{A}(t)^\top \\ \mathbf{A}(t) & \mathbf{0}_{q\times q}
		\end{bmatrix}.
		\end{align*}
		The strong convexity of $f_0(\bx,t)$ and the full-rank condition $\mbox{rank}(A(t))=q<n$ is sufficient for $\nabla_{\bz\bz} \mathcal{L}(\bz,t)$ to be invertible \cite{boyd2004convex}. Therefore, the Hessian inverse $\nabla_{\bz\bz}^{-1} {\mathcal{L}}(\bz,t)$ in \eqref{eq: time_varyin_equality_newton_method} exists. The time evolution of $\nabla_{\bz\bz} {\mathcal{L}}(\bz(t),t)$ can be written as	
		\begin{align*}
	    \dot{\nabla}_{\bz} \mathcal{L}(\bz(t),t)= \nabla_{\bz\bz} \mathcal{L}(\bz(t),t) \dot{\bz}(t) + \nabla_{\bz t} \mathcal{L}(\bz(t),t).
		\end{align*}
		Substituting $\dot{\bz}(t)$ from \eqref{eq: time_varyin_equality_newton_method}, we obtain
		\begin{align} \label{eq: Lagrangian_gradient_dynamics}
		\dot{\nabla}_{\bz} \mathcal{L}(\bz(t),t)=-\alpha \nabla_{\bz} \mathcal{L}(\bz(t),t),
		\end{align}
		which yields the solution
		\begin{align*}
		\nabla_{\bz} \mathcal{L}(\bz(t),t) = \nabla_{\bz} \mathcal{L}(\bz(0),0) e^{-\alpha t},
		\end{align*}
		for an initial condition $\bz(0) \in \mathbb{R}^{n+q}$. From the last identity, we obtain the bound
		\begin{align*}
		\|\nabla_{\bz} \mathcal{L}(\bz(t),t)\|_2 \leq \|\nabla_{\bz} \mathcal{L}(\bz(0),0)\|_2 e^{-\alpha t}.
		\end{align*}
		Hence, \eqref{eq: time_varying_equality_constraint_gradient_exponential} follows. Next, we apply the mean-value theorem to expand $\nabla_{\bz} \mathcal{L}(\bz(t),t)$ around the optimal point $\nabla_{\bz} \mathcal{L}(\bz^\star(t),t)=\mathbf{0}$ as follows,
			\begin{align}
			\nabla_{\bz} \mathcal{L}(\bz(t),t) =\nabla_{\bz\bz} \mathcal{L}(\boldsymbol{\eta}(t),t) (\bz(t)-\bz^\star(t)),\nonumber
			\end{align}
		where $\boldeta(t)$ is a convex combination of $\bz(t)$ and $\bz^\star(t)$. From the assumptions in the theorem, we have that $\|\nabla_{\bz\bz}^{-1} \mathcal{L}(\bz,t)\|_2 \leq M$ for all $z \in \mathbb{R}^{n+q}$ and $t \geq 0$. Therefore, we obtain from the last identity that
		\begin{align}
		\|\bz(t)-\bz^\star(t)\|_2 &= \|\nabla_{\bz\bz}^{-1} \mathcal{L}(\boldeta(t),t)\nabla_{\bz} \mathcal{L}(\bz(t),t)\|_2 \nonumber \\
		& \leq \|\nabla_{\bz\bz}^{-1} \mathcal{L}(\boldeta(t),t)\|_2 \|\nabla_{\bz} \mathcal{L}(\bz(t),t)\|_2 \nonumber \\
		& \leq M\|\nabla_{\bz} \mathcal{L}(\bz(0),0)\|_2 e^{-\alpha t}.\nonumber
		\end{align}
		On the other hand, recalling that $\bz(t)=[\bx(t)^\top \bnu^\top(t)]^\top$, we have the inequality $\|\bx(t)-\bx^\star(t)\|_2 \leq \|\bz(t)-\bz^\star(t)\|_2$. Combining these last two inequalities completes the proof.

		\subsection{Proof of Lemma \ref{lemma: pertubation_suboptimality_bounds}} \label{ap_perturbation_suboptimality_bound}
		Define $\widehat{\bbx}_s^{\star}(t)$ as 
		\begin{alignat}{2}
			\widehat{\bbx}_s^{\star}(t)
			:=\  &\underset{\bx \in \mathbb{R}^n}{\argmin}\ && f_{0}(\bbx,t)+\sum_{i=1}^{p} \mathbb{I}(s(t)-f_{i}(\bbx,t)),\nonumber
		\end{alignat}
		which is a perturbed version of the original optimization problem \eqref{eq: time_varying_problem_with_barrier} after including the slack variable $s(t)$ in the constraints. By perturbation and sensitivity analysis \cite[Chapter 5]{boyd2004convex}, we can establish the following inequality,
		\begin{align} \label{eq: slack_variable_perturbation}
		f_0(\widehat{\bbx}_s^{\star}(t),t) \leq f_0(\bx^{\star}(t),t) \leq f_0(\widehat{\bbx}_s^{\star}(t),t)+\sum_{i=1}^{p} \lambda_i^{\star}(t)s(t).
		\end{align} 				
		The first inequality is based on the fact that the feasible set is enlarged when $s(t) \geq 0$ and, hence, the optimal value is decreased. The second inequality follows directly from a sensitivity analysis of the original problem \cite[Chapter 5]{boyd2004convex}. On the other hand, replacing the indicator functions $\mathbb{I}(.)$ by a logarithmic barrier function, as in \eqref{eq: time_varying_problem_with_log_barrier_and_slack}, we obtain the bound \cite[Chapter 11]{boyd2004convex},
		\begin{equation} \label{eq: barrier_parameter_perturbation}
		f_0(\widehat{\bx}^{\star}(t),t) -  f_0(\widehat{\bx}_s^\star(t),t) \leq \frac{p}{c(t)}.
		\end{equation}
		It follows from \eqref{eq: slack_variable_perturbation}, \eqref{eq: barrier_parameter_perturbation}, and the triangle inequality that 
		$$
		|f_0(\widehat{\bbx}^{\star}(t),t)-f_0({\bx}^{\star}(t),t)| \leq \dfrac{p}{c(t)}+\sum_{i=1}^{p} \lambda^{\star}_i(t)s(t).
		$$
		The proof is complete.
		\subsection{Proof of Proposition \ref{lemma: time_varying_barrier_feasibility}} \label{ap_feasibility}
		For each constraint, we define the corresponding residual $\psi_i(\bx,t):=s(t)-f_i(\bx,t),\ i\in [p]$. Upon this definition, feasibility of $\bx(t)$ for all $t \geq 0$ is equivalent to non-negativity of the residuals at $(\bx(t),t)$ for all $t \geq 0$, i.e., $\psi_i(t)=\psi_i(\bx(t),t)\geq 0,\ \forall t\geq 0,\ i\in [p]$. We prove this by showing that 
		\begin{align} \label{eq: psi_limit}
		\lim\limits_{\psi_i \to 0^{+}} \frac{\dot{\psi}_i}{\psi_i}>0.
		\end{align}
		The above limit asserts that the time evolution of $\psi_i$ is strictly positive near the boundary. To prove \eqref{eq: psi_limit}, we first evaluate all the terms appearing in the ODE \eqref{eq: perturbed_time_varying_newton_barrier}. The gradient function $\nabla_{\bx}\Phi$ is
\begin{equation*}
\nabla_{\bx} {\Phi}=\nabla_{\bx} f_0+\dfrac{1}{c}\sum_{i=1}^{p} \dfrac{{\nabla_{\bx}f_i}}{\psi_i}.
\end{equation*}
The second partial derivatives read as
\begin{align*}
\nabla_{\bx\bx} {\Phi}&=\nabla_{\bx\bx} f_0+\dfrac{1}{c}\sum_{i=1}^{p} \left( \dfrac{\nabla_{\bx\bx}f_i}{\psi_i}
+\dfrac{\nabla_{\bx}f_i{\nabla_{\bx}f_i}^\top}{\psi_i^2}\right), \\
\nabla_{\bx s} \Phi &= \dfrac{1}{c}\sum_{i=1}^{p} -\dfrac{{\nabla_{\bx}f_i} }{\psi_i^2}, \\
\nabla_{\bx c} \Phi &= \dfrac{1}{c}\sum_{i=1}^{p} -\frac{1}{c}\dfrac{{\nabla_{\bx}f_i}}{\psi_i}, \\
\nabla_{\bx t} {\Phi}&=\nabla_{\bx t} f_0+\dfrac{1}{c}\sum_{i=1}^{p} \left( \dfrac{{\nabla_{\bx t}f_i}}{\psi_i} + \dfrac{{\nabla_{\bx}f_i} }{\psi_i^2}\frac{\partial f_i}{\partial t}\right).
\end{align*}
In the following, we study the limiting behavior of the above partial derivatives as $\psi_j \to 0$ for some $j \in [p]$. We have that

\begin{align} \label{eq: partial_derivatives_limit}
\psi_j \nabla_{\bx} {\Phi} &=  \dfrac{1}{c} \nabla_{\bx}f_j + O(\psi_j), \\
\psi_j \nabla_{\bx\bx} {\Phi}&=  \dfrac{1}{c} \dfrac{\nabla_{\bx}f_j{\nabla_{\bx}f_j}^\top}{\psi_j}+O(1), \nonumber\\
\psi_j \nabla_{\bx s} \Phi &= -\dfrac{1}{c}\dfrac{{\nabla_{\bx}f_j} }{\psi_j}+O(\psi_j), \nonumber\\
\psi_j \nabla_{\bx c} \Phi &= -\frac{1}{c^2}\nabla_{\bx}f_j + O(\psi_j), \nonumber\\
\psi_j \nabla_{\bx t} {\Phi}&=\dfrac{1}{c} \dfrac{{\nabla_{\bx}f_j}}{f_j} \frac{\partial f_j}{\partial t} + O(\psi_j). \nonumber
\end{align}
Multiply both sides of the ODE in \eqref{eq: perturbed_time_varying_newton_barrier} by $\psi_j$, and take the limit as $\psi_j \to 0$ to obtain
\begin{equation}\label{eqn_aux_lemma_5_proof}
\lim\limits_{\psi_j \to 0}\psi_j \Big[\nabla_{\bx \bx} \Phi \dot{\bx} +\alpha \nabla_{\bx} \Phi+ \nabla_{\bx s} \Phi \dot{s}+ \nabla_{\bx c}\dot{c}+\nabla_{\bx t} \Phi\Big]=0.
\end{equation}
Invoking \eqref{eq: partial_derivatives_limit} in the last identity yields,
\begin{align*}
\lim\limits_{\psi_j \to 0} \dfrac{\nabla_{\bx}f_j{\nabla_{\bx}^\top f_j}}{\psi_j} \dot{\bx}+\alpha \nabla_{\bx}f_j -\dfrac{ \nabla_{\bx}f_j }{\psi_j}\dot{s}+\dfrac{{\nabla_{\bx}f_j}}{f_j} \frac{\partial f_j}{\partial t} =0.
\end{align*}
Factoring out $\nabla_{\bx}f_j$, we obtain
\begin{align*}
\lim\limits_{\psi_j \to 0} \nabla_{\bx}f_j \Big[\dfrac{{\nabla_{\bx}f_j}^\top\dot{\bx}-\dot{s}+\dfrac{\partial f_j}{\partial t}}{\psi_j}+\alpha\Big]=0.
\end{align*}
\ifx
Using the expressions derived in \eqref{eqn_aux_lemma_5_proof_grad}, \eqref{eqn_aux_lemma_5_proof_gradt} and \eqref{eqn_aux_lemma_5_proof_hess}, and simplifying the resulting expression the limit of \eqref{eqn_aux_lemma_5_proof} when $\psi_j \to 0$ reduces to
\begin{equation*}
\lim_{\psi_j\to0} \nabla_{\bx}f_j(\bbx,t)  \dfrac{\frac{\partial \psi_j}{\partial t} (\bbx,t)-{\nabla_{\bx}f_j}^\top(\bbx,t) \dot{\bx}}{\psi_j(\bbx,t)}=\alpha\nabla_{\bx}f_j(\bbx,t). 
\end{equation*} 
\fi
Using the identity $\dot{\psi_i}=\dot{s}-\dfrac{\partial f_j}{\partial t} -{\nabla_{\bx}f_j}^\top \dot{\bx}$, we conclude that
		\begin{align*}
		\lim\limits_{\psi_j \to 0^+}\dfrac{\dot{\psi}_j}{\psi_j} = \alpha.
		\end{align*}
Since $\alpha>0$ we showed that $\psi_j$ is strictly increasing at the boundary of the feasible set, hence ensuring non-negativity of the residuals $\psi_i(\bbx(t),t)$ for all $t \geq 0$ and $i \in [p]$. As explained at the beginning of this proof, non-negativity of the residuals is equivalent to feasibility of $\bbx(t)$ for all $t \geq 0$. 		
		
		\subsection{Proof of Lemma \ref{lemma: time_varying_barrier_convergence}} \label{ap: lemma_time_varying_barrier_convergence}
				Since $f_0(\bx,t)$ is strongly convex, and $c(t)$ is strictly positive, it follows that $\nabla_{\bx\bx}{\Phi}$ is $m$-strongly convex for $\bx \in \mathcal{D}(t)$ and, therefore, $\nabla_{\bx\bx}^{-1} \Phi$ exists and is bounded. The dynamics of $\nabla_{\bx} {\Phi}$ can be written as $\dot\nabla_{\bx} {\Phi} =\nabla_{\bx\bx} {\Phi}\dot\bx+ \nabla_{\bx t} {\Phi}$.
				Substituting $\dot{\bx}$ from \eqref{eq: perturbed_time_varying_newton_barrier} into the last result results in the closed-loop dynamics $\dot\nabla_{\bx} {\Phi} =-\alpha\nabla_{\bx} {\Phi} $, which in turn implies that
				\begin{equation*}
				\|\nabla_{\bx} {\Phi}({\bx}(t),s(t),c(t),t)\|_{2} \leq e^{-\alpha t}\|\nabla_{\bx} {\Phi}(\bx(0),s(0),c(0),0)\|_2.
				\end{equation*}
				Finally, it follows from $m$-strong convexity of $\Phi$ that
				\begin{align*} 
				\|{\bx}(t)-\widehat{\bx}^{\star}(t)\|_2 \leq \dfrac{1}{m} \|\nabla_{\bx} {\Phi}({\bx}(t),s(t),c(t),t)\|_{2}.
				\end{align*}
				Combining the last two inequalities gives the desired inequality. The proof is complete.
				
				\ifx
				\subsection{Proof of Theorem \ref{thm_main_theorem}} \label{ap: thm_main_theorem}
					By Assumption \ref{assumption: bounded_in_times}, we have that $\lim\limits_{t \to \infty}  \sum_{i=1}^{p} \lambda_i^\star(t) s(t)$ when $s(t)=s(0)\exp(-\gamma_s t)$. This together with $\lim\limits_{t \to \infty} c(t)=\infty$ implies that $\lim\limits_{t \to \infty} |f_0(\bx^\star(t),t) - f_0(\hat{\bx}^\star_s(t),t)|=0$. Strong convexity of $f_0$ would let us to conclude that $\lim\limits_{t \to \infty }\|\widehat{\bx}^\star_s(t)-\bx^\star(t)\|_2 \to 0$. The proof is complete.
				\fi		
					\subsection{Proof of Theorem \ref{lem: unconstrained_finite_time}} \label{ap: lem: unconstrained_finite_time}
						We first define the following Lyapunov function,
						\begin{align}
						V(\bx,s,c,t) = \frac{1}{2} \|\nabla_{\bx} \Phi(\bx,s,c,t)\|_2^2,
						\end{align}
						which is positive everywhere and is zero along the approximate optimal trajectory, i.e., $V(\widehat{\bx}^\star(t),s(t),c(t),t)=0,\ t\geq 0$. The time derivative of the Lyapunov function along the trajectories of \eqref{eq: perturbed_time_varying_newton_barrier_robust} is
						\begin{align}
						\dot{V} &= \nabla_{\bx} \Phi^\top \dot{\nabla}_{\bx} \Phi \nonumber\\
						&=\nabla_{\bx} \Phi^\top (\nabla_{\bx\bx} \Phi\, \dot{\bx}+\nabla_{\bx s} \Phi \dot{s}+\nabla_{\bx c} \Phi \dot{c} + \nabla_{\bx t} \Phi) \nonumber \\
						&=\nabla_{\bx} \Phi^\top(-\alpha\nabla_{\bx} \Phi+\nabla_{\bx t} \Phi-\widehat{\nabla}_{\bx t} \Phi).\nonumber
						\end{align}
						When $\|\nabla_{\bx} \Phi\|_2 \geq \varepsilon$, we have that $\alpha=\frac{\alpha_0}{\|\nabla_{\bx} \Phi\|_2}$ and, therefore, $\dot{V}$ is given by
						\begin{align*}
						\dot{V}&=\nabla_{\bx} \Phi^\top (-\alpha \dfrac{\nabla_{\bx} \Phi}{\|\nabla_{\bx} \Phi\|_2}-\widehat\nabla_{\bx t} \Phi+\nabla_{\bx t} \Phi)\\
						&=-\alpha_0 \|\nabla_{\bx} \Phi\|_2+\nabla_{\bx} \Phi^\top (-\widehat\nabla_{\bx t} \Phi+\nabla_{\bx t} \Phi).
						\end{align*}
						Using the Assumption $\|\widehat\nabla_{\bx t} \Phi-\nabla_{\bx t} \Phi\|_2 \leq \eta$, we obtain the inequality
						\begin{align}
						\dot{V} &\leq (\eta-\alpha_0) \|\nabla_{\bx} \Phi\|_2 \nonumber\\
						&= (\eta-\alpha_0) \sqrt{2V}.
						\end{align}
						Using the comparison lemma \cite{khalil1996nonlinear}, we can write $V(t) \leq W(t)$, where $W(t)$ is the solution of the initial value problem $\dot W(t)=(\eta-\alpha_0)\sqrt{2W(t)},\ W(0)=V(0)$. From the last ODE, we obtain the solution $W(t)=\frac{1}{2}(\sqrt{2V(0)}-(\alpha_0-\eta)t)^2$. Hence, the Lyapunov function satisfies the bound
						\begin{align*}
						2V(t) \leq (\sqrt{2V(0)}-(\alpha_0-\eta)t)^2,
						\end{align*} 
						or, equivalently,
						\begin{align}
						\|\nabla_{\bx}\Phi\|_2 \leq  \|\nabla_{\bx} \Phi_0 \|_2 - (\alpha_0-\eta) t.\nonumber
						\end{align}
						The right-hand side becomes equal to $\varepsilon$ in finite time, implying that the trajectory reaches the set $\mathcal{S}_t(\varepsilon)$ in finite time. When $\|\nabla_{\bx} \Phi\|_2 \leq \varepsilon$, we have that $\alpha=\alpha_0/\varepsilon$, and the time derivative of the Lyapunov function becomes
						\begin{align*}
						\dot{V}&=\nabla_{\bx} \Phi^\top (-\frac{\alpha_0}{\varepsilon} \nabla_{\bx} \Phi-\widehat\nabla_{\bx t} \Phi+\nabla_{\bx t} \Phi)\\
						&=-\frac{\alpha_0}{\varepsilon} \|\nabla_{\bx} \Phi\|_2^2+\nabla_{\bx} \Phi (-\widehat\nabla_{\bx t} \Phi+\nabla_{\bx t} \Phi)\\
						& \leq \|\nabla_{\bx} \Phi\|_2 (-\frac{\alpha_0}{\varepsilon} \|\nabla_{\bx} \Phi\|_2+\eta).
						\end{align*}	
						It is evident that $\dot{V}$ is negative when $\|\nabla_{\bx} \Phi\|_2>\eta \varepsilon/\alpha_0$, or equivalently, $\dot{V}$ is negative outside the set $S_t(\eta \varepsilon/\alpha_0)$. As a result, the solution $\bx(t)$ converges asymptotically to this set. Notice that since $\eta<\alpha_0$, we have that $S_t(\eta  \varepsilon/\alpha_0) \subset S_t(\varepsilon)$. In other words, the solution converges to $\mathcal{S}_t(\varepsilon)$ in finite time, and stays there forever. The proof is complete.	
				
\subsection{Expressions for the Numerical Examples }\label{ap_numerical_examples}
We now derive explicit expressions for the terms in the ODE \eqref{eq: perturbed_time_varying_newton_barrier} for the example in Subsection \ref{sec_robot_navigation}. The gradient of the augmented objective function with respect to $\bbx$ takes the form
\begin{equation}
\nabla_{\bbx}\Phi(\bbx,\bbx_c,t) = \bbx-\bbx_d +\frac{1}{c(t)}\sum_{i=1}^m\frac{\bba_i(\bbx_c)}{b_i(\bbx_c)-\bba_i(\bbx_c)^\top\bbx},\nonumber
\end{equation}
and its Hessian reads as
\begin{equation}
\nabla_{\bbx \bbx}\phi(\bbx,\bbx_c,t)(\bbx,\bbx_c) = \bbI_n +\frac{1}{c(t)}\sum_{i=1}^m\frac{\bba_i(\bbx_c)\bba_i(\bbx_c)^\top}{(b_i(\bbx_c)-\bba_i(\bbx_c)^\top\bbx)^2}.\nonumber
\end{equation}
Furthermore, the time derivative of the gradient of the barrier function can be written as
\begin{equation}
\begin{split}
\nabla_{\bbx t}\phi(\bbx,\bbx_c,t) = -\frac{\dot{c}(t)}{c(t)^2}\sum_{i=1}^m\frac{\bba_i(\bbx_c)}{b_i(\bbx_c)-\bba_i(\bbx_c)^\top\bbx} \\
+\frac{1}{c(t)}\sum_{i=1}^m\frac{\dot{\bba}_i(\bbx_c)}{b_i(\bbx_c)-\bba_i(\bbx_c)^\top\bbx} \nonumber\\
-\frac{1}{c(t)}\sum_{i=1}^m\bba_i(\bbx_c)\frac{\dot{b}_i(\bbx_c)-\dot{\bba}_i(\bbx_c)^\top\bbx}{(b_i(\bbx_c)-\bba_i(\bbx_c)^\top\bbx)^2}.
\end{split}
\end{equation}
The expressions for $\dot{\bba}_i(\bbx_c)$ and $\dot{b}_i(\bbx_c)$ are derived below. The expression of $\bba_i(\bbx_c)$ for every $i=[m]$ is given by $\bba_i(\bbx_c) = (\bbx_i - \bbx_c)$. Thus, its time derivative is given by 
\begin{equation}
\dot{\bba}_i(\bbx_c) = - \dot{\bbx}_c = K(\bbx_c-\hat{\bbx}),\nonumber
\end{equation}
where the last equality comes from replacing the time derivative of $\bbx_c$ by the control law \eqref{eqn_control_law}. We derive next the expression for the time derivative of $b_i(\bbx_c)$ defined in \eqref{eqn_bi}.
To do so, we compute the time derivative of $\theta_i(\bbx_c)$. Differentiating $\theta_i(\bbx_c)$, defined in \eqref{eqn_bi}, yields
\begin{equation}\label{eqn_alphadot}
\dot{\theta}_i(\bbx_c) = -\frac{(\bbx_c-\bbx_i)^\top\dot{\bbx}_c}{\|\bbx_c-\bbx_i\|^4} = K\frac{(\bbx_c-\bbx_i)^\top(\bbx_c-\hat{\bbx})}{\|\bbx_c-\bbx_i\|^4}.
\end{equation}
Differentiating $b_i(\bbx)$ in \eqref{eqn_bi} yields 
\begin{align}
\dot{b}_i(\bbx_c) =& -\dot{\bbx}_c^\top\left( \theta_i\bbx_i+(1-\theta_i)\bbx_c \right) \nonumber\\
 &+\dot{\theta}_i\|\bbx_i-\bbx_c\|^2(r-\theta_i)(\bbx_i-\bbx_c)^\top\dot{\bbx}_c,\nonumber
\end{align}
where in the above equation, $\dot{\theta}_i$ and $\dot{\bbx}_c$ are respectively given by \eqref{eqn_alphadot} and \eqref{eqn_control_law}.

\bibliographystyle{ieeetr}
\bibliography{Refs}

\begin{IEEEbiography}[{\includegraphics[width=1in,height=1.25in,clip,keepaspectratio]{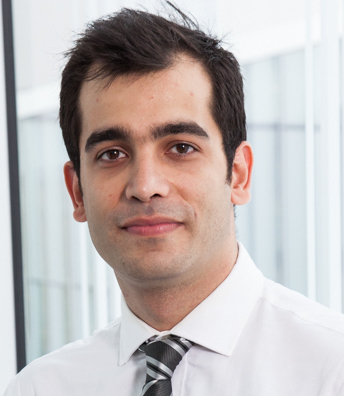}}]{Mahyar Fazlyab}
	received his B.Sc. and M.Sc. degrees in Mechanical Engineering from Sharif University of Technology, Tehran, Iran, in 2010 and 2013. He has been a PhD student with the Department of Electrical and Systems Engineering at the University of Pennsylvania since September 2013. His research interests include the analysis, optimization, and control of (networked) dynamical systems.
\end{IEEEbiography}

\begin{IEEEbiography}[{\includegraphics[width=1in,height=1.25in,clip,keepaspectratio]{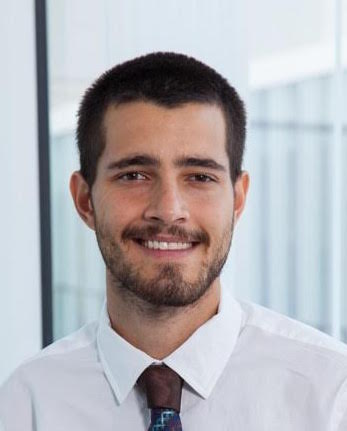}}]{Santiago Paternain}
	received the B.Sc. degree in Electrical Engineering from Universidad de la Rep\'ublica Oriental del Uruguay, Montevideo, Uruguay in 2012. Since August 2013, he has been working toward the Ph.D. degree in the Department of Electrical and Systems Engineering, University of Pennsylvania. His research interests include optimization and control of dynamical systems. 
\end{IEEEbiography}

\begin{IEEEbiography}[{\includegraphics[width=1in,height=1.25in,clip,keepaspectratio]{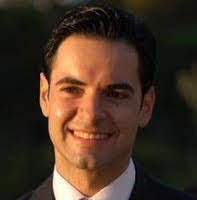}}]{Victor M.~Preciado} 
	received his Ph.D. degree in Electrical Engineering and Computer Science from the Massachusetts Institute of Technology in 2008. He is currently the Raj and Neera Singh Assistant Professor of Electrical and Systems Engineering at the University of Pennsylvania. He is a member of the Networked and Social Systems Engineering (NETS) program and the Warren Center for Network and Data Sciences. His research interests include network science, dynamic systems, control theory, and convex optimization with applications in socio-technical systems, technological infrastructure, and biological networks.
\end{IEEEbiography}

\begin{IEEEbiography}[{\includegraphics[width=1in,height=1.25in,clip,keepaspectratio]{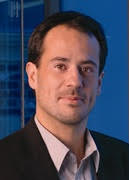}}]{Alejandro Ribeiro}
	received the B.Sc. degree in electrical engineering from the Universidad de la Republica Oriental del Uruguay, Montevideo, in 1998 and the M.Sc. and Ph.D. degree in electrical engineering from the Department of Electrical and Computer Engineering, the University of Minnesota, Minneapolis in 2005 and 2007. From 1998 to 2003, he was a member of the technical staff at Bellsouth Montevideo. After his M.Sc. and Ph.D studies, in 2008 he joined the University of Pennsylvania (Penn), Philadelphia, where he is currently the Rosenbluth Associate Professor at the Department of Electrical and Systems Engineering. His research interests are in the applications of statistical signal processing to the study of networks and networked phenomena. His focus is on structured representations of networked data structures, graph signal processing, network optimization, robot teams, and networked control. Dr. Ribeiro received the 2014 O. Hugo Schuck best paper award, the 2012 S. Reid Warren, Jr. Award presented by Penn's undergraduate student body for outstanding teaching, the NSF CAREER Award in 2010, and paper awards at the 2016 SSP Workshop, 2016 SAM Workshop, 2015 Asilomar SSC Conference, ACC 2013, ICASSP 2006, and ICASSP 2005. Dr. Ribeiro is a Fulbright scholar and a Penn Fellow.
\end{IEEEbiography}

\end{document}